\pdfoutput=1

\documentclass[a4paper, 11pt, abstracton, titlepage, oneside]{scrartcl}
\makeatletter

\usepackage[T1]{fontenc}
\usepackage[utf8]{inputenc}
\usepackage[onehalfspacing]{setspace}
\usepackage[headsepline]{scrlayer-scrpage}
\clearscrheadfoot 
\usepackage{layout}
\usepackage{geometry}
\usepackage{adjustbox}
\usepackage{authblk}
\usepackage[titletoc,title]{appendix}

\usepackage[hyphens]{url}
\usepackage{color}

\usepackage{amsmath,amssymb,amsfonts,amsthm, mathtools, bm}
\usepackage{thmtools, thm-restate}
\usepackage{nicefrac}
\usepackage{array}

\usepackage[normal]{threeparttable}
\usepackage{threeparttablex}
\usepackage{tabularx}
\usepackage{longtable}
\usepackage{booktabs}
\usepackage{colortbl}
\usepackage{multirow}
\usepackage{makecell}

\usepackage[acronym]{glossaries}

\usepackage{subcaption}
\usepackage{floatrow}
\usepackage{graphicx}
\usepackage{placeins}

\usepackage{tikz}
\usetikzlibrary{arrows.meta}
\usetikzlibrary{math}
\usetikzlibrary{shapes}
\usepackage{pgfplots}
\usepgfplotslibrary{groupplots}
\usepgfplotslibrary{dateplot}
\usepackage{tikzscale}

\usepackage{enumerate}
\usepackage{multicol}
\usepackage[author=Patrick]{pdfcomment}
\usepackage{xparse}
\usepackage{twoopt}
\usepackage{etoolbox}

\usepackage{eurosym}
\usepackage{ellipsis}
\usepackage{natbib}
\usepackage{paralist}
\usepackage{ulem}

\usepackage[ruled,linesnumbered]{algorithm2e}
\usepackage{bbm}
\AtBeginDocument{
	\newgeometry{
		textheight=9in,
		textwidth=6.5in,
		top=1in,
		headheight=14pt,
		headsep=25pt,
		footskip=30pt
	}
}
\newsavebox{\abstractbox}
\renewenvironment{abstract}
{\begin{lrbox}{0}\begin{minipage}{\textwidth}
			\begin{center}\normalfont\sectfont\abstractname\end{center}\quotation}
		{\endquotation\end{minipage}\end{lrbox}%
	\global\setbox\abstractbox=\box0 }

\makeatletter
\expandafter\patchcmd\csname\string\maketitle\endcsname
{\vskip\z@\@plus3fill}
{\vskip\z@\@plus2fill\box\abstractbox\vskip\z@\@plus1fill}
{}{}
\makeatother

\addtokomafont{captionlabel}{\footnotesize\bfseries} 
\addtokomafont{caption}{\footnotesize\bfseries} 

\bibpunct[, ]{(}{)}{,}{a}{}{,}%
\newtheorem{definition}{Definition}[section]
\newtheorem{theorem}{Theorem}[section]

\let\endproof\relax

\counterwithin*{equation}{section}

\newenvironment{myfont}{\fontfamily{ccr}\selectfont}{\par}
\DeclareTextFontCommand{\textmyfont}{\myfont}

\AtBeginDocument{%
	\abovedisplayskip=7.5pt plus 0pt minus 0pt
	\abovedisplayshortskip=7.5pt plus 0pt minus 0pt
	\belowdisplayskip=7.5pt plus 0pt minus 0pt
	\belowdisplayshortskip=7.5pt plus 0pt minus 0pt
}

\newcolumntype{L}[1]{>{\raggedright\let\newline\\\arraybackslash\hspace{0pt}}p{#1}}
\newcolumntype{C}[1]{>{\centering\let\newline\\\arraybackslash\hspace{0pt}}p{#1}}
\newcolumntype{R}[1]{>{\raggedleft\let\newline\\\arraybackslash\hspace{0pt}}p{#1}}
\renewcommand{\emph}[1]{\textit{#1}}
\graphicspath{{Figures/}}

\begin{document}
\emergencystretch 3em
\title{\large Branch-and-Price for the Stochastic TSP with Generalized Latency}

\author[1]{\normalsize Benedikt Lienkamp}
\author[2]{\normalsize Mike Hewitt}
\author[3]{\normalsize Maximilian Schiffer}
\affil{\small 
	TUM School of Management, Technical University of Munich, 80333 Munich, Germany
	
	\scriptsize benedikt.lienkamp@tum.de

    \small
    \textsuperscript{2}Quinlan School of Business, Loyola University Chicago, USA

    \scriptsize mhewitt3@luc.edu
    
	\small
	\textsuperscript{3}TUM School of Management \& Munich Data Science Institute,
	
	Technical University of Munich, 80333 Munich, Germany
	
	\scriptsize schiffer@tum.de}

\date{}

\lehead{\pagemark}
\rohead{\pagemark}

\begin{abstract}
\begin{singlespace}
{\small\noindent Motivated by the tactical planning level of demand adaptive public transportation systems, we present the stochastic symmetric traveling salesman problem with generalized latency (STSP-GL), a stochastic extension to the symmetric traveling salesman problem with generalized latency (TSP-GL). The STSP-GL aims to choose a subset of nodes of an undirected graph and determines a Hamiltonian tour amongst those nodes, minimizing an objective function that is a weighted combination of route design and passenger routing costs. These nodes are selected to ensure that a predefined percentage of uncertain passenger demand is served with a given probability. We formulate the STSP-GL as a stochastic program and propose a branch-and-price algorithm for solving its deterministic equivalent. We also develop a local search approach with which we improve the performance of the B\&P approach. We assess the efficiency of the proposed methods on a set of instances from the literature. We demonstrate that the proposed methods outperform a known benchmark, improving upper bounds by up to 85\% and lower bounds by up to 55\%. Finally, we show that solutions of the stochastic model are both more cost-effective and robust than those of the deterministic model.
\\
\medskip}
\end{singlespace}
\end{abstract}

\maketitle
\setcounter{equation}{0}
\section{Introduction}
As urban areas experience rapid population growth, the importance of public transportation becomes increasingly evident. In light of expanding cities, challenges such as traffic congestion, environmental pollution, and limited resources become more pressing. Within this context, public transportation takes on a crucial role as it offers a viable solution to tackle these challenges by presenting an alternative to private vehicle usage. This alternative helps to alleviate traffic congestion and to reduce pollution. Moreover, public transportation ensures accessibility for individuals across diverse socioeconomic backgrounds, connecting them to vital services, education, and employment opportunities. Additionally, its ability to optimize land use and curb urban sprawl promotes compact, well-connected communities while enhancing social cohesion. By establishing and strengthening comprehensive public transportation networks, cities can attain improved organizational efficiency, heightened environmental sustainability, and greater societal inclusivity.

Based on these considerations, traditional fixed bus lines serve as a cost-effective and reliable foundation of public transportation networks, offering widespread coverage that connects diverse communities. Their high passenger capacity and established routes promote social inclusivity and accessibility to essential services. Moreover, they enable efficient utilization of existing infrastructure, allowing to transport a larger number of people using fewer vehicles, which can contribute to reduced traffic congestion and emissions. However, they can be susceptible to traffic congestion and limited flexibility in adapting to changing travel demands, leading to potential delays and less personalized service experiences for passengers.

Demand Responsive Systems (DRS) offer more flexibility and can adapt to fluctuating individual transportation demands, thus enabling more personalized services while maintaining a degree of resource sharing. DRS were first introduced under the name of Dial-a-Ride (DAR) as a door-to-door service for users with reduced mobility (\citealp{wilson1971scheduling}; \citealp{ioachim1995request}; \citealp{toth1996fast}). One drawback of DAR are the higher operational costs associated with providing individualized service and accommodating specific travel requests. Additionally, the real-time nature of pickups and drop-offs, i.e., significant route changes during operation, can lead to unpredictable travel times for passengers, which might not be suitable for those with time-sensitive commitments.

\begin{figure}[t]
\centering
\begin{subfigure}{0.45\textwidth}
  \centering
    \includegraphics[width=\columnwidth]{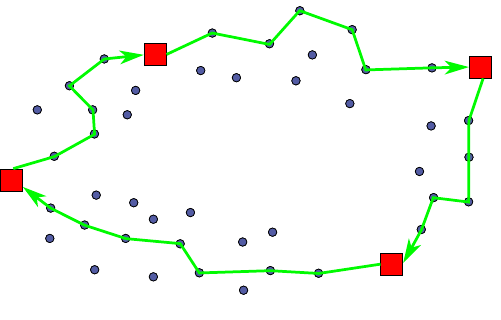}
    \captionsetup{format=hang}
	\caption{Route variant 1}
\end{subfigure}
\begin{subfigure}{0.45\textwidth}
  \centering
    \includegraphics[width=\columnwidth]{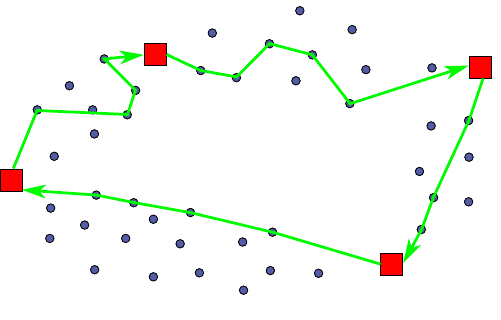}
    \captionsetup{format=hang}
    \caption{Route variant 2}
\end{subfigure}
\caption{A DAS line with compulsory stops (red) and optional stops (blue). Picture modified from \cite{crainic2012designing}}
\label{fig: DAS line}
\end{figure}

Demand Adaptive Systems (DAS) (\citealp{malucelli1999demand}; \citealp{quadrifoglio2007insertion}) incorporate DRS into scheduled bus transportation.
To do so, a DAS bus line provides a traditional transit-line service for a set of compulsory stops. These compulsory stops are bound to a schedule with fixed time windows during which the vehicle serving the line has to leave each compulsory stop. Additionally, passengers may issue requests at stops that are not compulsory, inducing detours in the vehicle routes (see Figure~\ref{fig: DAS line}). We refer to these non-compulsory stops as optional stops. The underlying principle of DAS is that the regularity and predictability of the service is a critical attribute of public transportation, as it enables passengers to plan their trips, facilitates integration with other modes of transportation, and allows access to the service without the need for advanced booking. This regularity is considered a valuable property of public transit as it enhances the user experience and enables more efficient use of the transportation infrastructure. Additionally, the flexibility of DAS allows for a more personalized transportation system with the possibility of serving a higher share of passengers than traditional fixed bus lines at a lower cost than DAR. This combination of flexibility and predictability merges the advantages of traditional fixed bus lines and DAR while limiting their drawbacks.  

Planning the operations of DAS is a highly complex task as it requires balancing the characteristics of traditional and on-demand systems. This involves both a service design phase, where routes and schedules are established, and an operational phase, with adjustments to vehicle routes and schedules in response to user requests. The complexity in planning DAS operations has been extensively studied in the literature, with a comprehensive review of methodological aspects presented in \cite{errico2013survey}. Focusing on the aforementioned service design phase, the tactical planning of a DAS line encompasses three main interrelated activities: (1) identification and selection of transit stops to be visited in accordance with a pre-determined, regular timetable (compulsory stops); (2) determination of a general order of transit stops within the service area; and (3) establishment of schedules for compulsory stops. In this paper, we focus on activity (2), and, unlike existing literature, we explicitly model uncertainty in passenger demands, which allows us to determine a general order of transit stops that is robust against volatile passenger demand.

Given that the compulsory stops have been identified, activity (2) can be modeled as a symmetric traveling salesman problem with generalized latency. This problem has been studied in \cite{errico2017benders} and aims to find a route between all compulsory and optional transit stops, simultaneously minimizing passenger travel times and company operating costs. \cite{errico2011design} describe the emergence and interaction of the TSP-GL within the broader context of DAS tactical planning. Through experimental analysis, they emphasize that passenger perceptions of service quality depend greatly on travel times.
\cite{errico2017benders} propose a methodological approach for the optimization of the TSP-GL through the use of a Benders-branch-and-cut (BBC) algorithm. In this work, the authors assumed deterministic passenger demand. Building on this methodological foundation, we aim to integrate stochastic passenger demand into this algorithmic framework to account for the uncertainty of passenger demand during the tactical planning of DAS and optimize the resulting stochastic TSP-GL (STSP-GL). Additionally, we do not assume all demand must be met but instead seek to meet a service level. These extensions allow DAS operators to analyze the impact of different service levels and the frequency of meeting them on the route design and served passenger demand. Furthermore, accounting for uncertainty in passenger demand allows for the design of DAS routes that are robust against volatile passenger demand.
\subsection{Literature Review}
We define the STSP-GL as a generalization of the TSP-GL in which only a subset of all stochastic passenger demand can be served as long as a service level that measures the expected sum of all transported demand is met with a certain probability. Like the TSP-GL, the STSP-GL is closely related to the Fixed-Charge Network Design Problem (FNDP), which involves finding a set of arcs in a graph that enables the flow of commodities from their origins to their destinations in order to satisfy some demand characteristics. The FNDP minimizes routing costs, i.e., variable expenses incurred when sending flows through the network, or design costs, i.e., fixed expenses associated with establishing network infrastructure, or both \citep{costa2005survey}. In the STSP-GL, all arcs of the graph have infinite capacity, and the set of chosen arcs has to form a tour that visits all compulsory stops. Additionally, this set of arcs has to allow passenger flow for a subset of passengers to meet a predefined service level, i.e., the expected sum of served passenger demand has to be met with a certain probability. Several exact and heuristic methods have been proposed to solve the FNDP.
The most prominent exact methods to solve FNDP are Lagrangian-based techniques, Benders decomposition, and polyhedral approaches. For an extensive overview of exact methods for FNDP, we refer to \cite{gendron1999multicommodity}, \cite{gendron2011decomposition}, and \cite{crainic2021fixed}. 

Lagrangian-based techniques utilize either subgradient algorithms or bundle methods to solve the Lagrangian dual. 
\cite{crainic2001bundle} show that for Lagrangian-based techniques, bundle methods appear superior to subgradient approaches as they converge faster and are more robust relative to different relaxations and problem characteristics. \cite{frangioni2014bundle} exploit the structure of the Dantzig-Wolfe reformulation of the Lagrangian dual to construct an efficient bundle method. \cite{frangioni2017computational} analyzed a large class of subgradient methods to show that the theory developed in \cite{d2009convergence} impacts the performance of subgradient methods. \cite{gendron2019revisiting} compares Lagrangian relaxation bounds for a large class of network design problems, which they obtain either by relaxing the linking constraints or the flow conservation constraints.

\cite{costa2009benders} give an extensive overview on the application of Benders decomposition to FNDP and clarify the relationships between Benders, metric, and cutset inequalities, which are often used in solution algorithms for multi-commodity capacitated FNDP. \cite{costa2012accelerating} propose a general scheme for generating extra cuts in a Benders decomposition, which they apply to the FNDP. 

There exists a multitude of polyhedral approaches to solve the FNDP. \cite{wolsey2003strong} gives an overview on strong valid inequalities and tight extended formulations for general mixed integer linear programs. \cite{chouman2017commodity} develop a cutting-plane algorithm for the multi-commodity capacitated FNDP. \cite{atamturk2002capacitated} analyzes the cut-set polyhedra of capacitated FNDP and shows that utilizing cut-set inequalities significantly improves computations. \cite{atamturk2001flow} applies flow pack inequalities, which improve computations for capacitated fixed-charge network flow problems. \cite{padberg1985valid} study the single-node fixed-charge flow problem and introduce flow cover inequalities.

The most prominent heuristic approaches to solve FNDP include classical heuristics such as local search \citep{walker1976heuristic}, neighborhood-based metaheuristics such as tabu search \citep{Gendreau2019}, simulated annealing \citep{delahaye2019simulated}, and iterated local search \citep{lourencco2019iterated}, population-based metaheuristics such as genetic algorithms \citep{whitley2019next}, and parallel metaheuristics \citep{crainic2019parallel}. Furthermore, several matheuristics have been applied to FNDP. 

There exist several approaches to solve FNDP via tabu search. \cite{sun1998tabu} apply tabu search with two strategies for each of the intermediate and long-term memory processes to the FNDP.
\cite{crainic2000simplex} utilize a tabu search framework that combines pivot moves with column generation to solve multicommodity capacitated FNDP. \cite{ghamlouche2003cycle} propose cycle-based neighborhood structures for multicommodity-capacitated FNDP and show their efficiency through a simple tabu search method.

For neighborhood-based metaheuristics, \cite{lotfi2013genetic} introduce a priority-based encoding genetic algorithm for linear and nonlinear FNDP and show that it outperforms a spanning tree-based genetic algorithm in terms of solution quality and computation time. \cite{ghamlouche2004path} combine path relinking with cycle-based tabu search in an algorithm that is robust regarding solution quality and computational effort. \cite{doerner2007} propose a scatter search approach for the capacitated FNDP, which cannot match the path relinking algorithm consistently but can produce better solutions on some of the test instances.

Several matheuristics have been applied to FNDP. \cite{rodriguez2010local} utilize local branching \citep{fischetti2003local} to solve the multicommodity capacitated FNDP. They show that their algorithm outperforms the cycle-based tabu search \citep{ghamlouche2003cycle} and path relinking \citep{ghamlouche2004path}. 

\cite{hewitt2010combining} combine mathematical programming algorithms and heuristic search techniques to find provably high-quality solutions for capacitated FNDP quickly. They show that their algorithm outperforms cycle-based tabu search \citep{ghamlouche2003cycle} and path relinking \citep{ghamlouche2004path} as well.

\cite{katayama2011combining} combine capacity scaling using a column-row generation technique and local branching. They show that their approach outperforms cycle-based tabu search \citep{ghamlouche2003cycle}, path relinking \citep{ghamlouche2004path}, local branching \citep{rodriguez2010local}, and IP Search \citep{hewitt2010combining}.

\cite{yaghini2015cutting} propose a cutting-plane neighborhood structure, which they implement in a tabu search framework and combine with local branching. They show that their approach outperforms cycle-based tabu search \citep{ghamlouche2003cycle}, path relinking \citep{ghamlouche2004path}, local branching \citep{rodriguez2010local}, IP Search \citep{hewitt2010combining} and CALB \citep{katayama2011combining}. 
We refer to \cite{crainic2021fixed} for an extensive overview of all exact and heuristic approaches to solve the FNDP.

Besides FNDP, the STSP-GL relates to three additional problems, the generalized minimum latency problem (GMLP), the TSP-GL, and the selective general minimum latency problem (SGMLP). \cite{errico2008design} introduces the GMLP and the SGMLP. The GMLP considers a complete directed graph $G = (N, A)$ with associated design and routing cost on its directed arcs as well as deterministic origin-destination demands $d_{hk}$ between its node pairs. The goal is to identify a Hamiltonian directed cycle that minimizes a combination of design costs, i.e., the time needed to travel the tour, and routing costs, i.e., average passenger travel time, of the demand on this cycle. The SGMLP builds on this setting and relaxes the condition of the Hamiltonian tour. Accordingly, the goal of the SGMLP is to find a simple tour on the graph, which does not have to be Hamiltonian, i.e., to serve a subset of deterministic passenger demands such that a given service level is met. The TSP-GL also builds on the methodology of the GMLP with the difference in the cycle design, as the GMLP seeks a directed tour while the TSP-GL seeks an undirected tour (\citealp{errico2008design}; \citealp{errico2017benders}). This difference allows the GMLP to consider asymmetric design costs but also fixes the flow direction of the passengers. The STSP-GL combines the TSP-GL and the SGMLP as it seeks an undirected tour that does not have to include all graph nodes. Furthermore, it also assumes stochastic passenger demand instead of deterministic passenger demand, which is the case for GMLP, SGMLP, and TSP-GL. \cite{errico2017benders} give an overview on how the TSP-GL and, thus also, the STSP-GL additionally relate to the minimum latency problem.

The GMLP, SGMLP, and TSP-GL have all been studied under the assumption that all information about passenger demand is available before the design decisions, i.e., which edges form the tour, are made. It is a general understanding that, in most cases, this does not translate to real-world applications \citep{kall1994stochastic}. Traditionally, this is often not considered explicitly during the tactical planning phase but dealt with in the operational planning phase \citep{ibarra2015planning}. Accordingly, passenger demand is usually determined through forecasting methods, e.g., average historical data. Usually, the expected quality of a solution to a stochastic model is better than a solution to the deterministic counterpart when evaluated in the stochastic setting. \cite{wallace2000decision} and \cite{higle2003sensitivity} explain that this occurs, as the deterministic model is optimal for one specific scenario but might be very bad for scenarios where it is not optimal. Additionally, the deterministic solution is often feasible but not optimal in the stochastic model. \cite{birge1982value} measures this difference in solution quality by the value of stochastic solution (VSS). Here, the VSS represents the potential benefit of solving the stochastic model instead of the deterministic model in the stochastic environment. \cite{wallace2000decision} shows that stochastic models may find solutions significantly different from deterministic models. Nevertheless, \cite{thapalia2011stochastic} and \cite{maggioni2012analyzing} show that for multiple other problems, certain structural patterns from the deterministic solutions may re-emerge in the stochastic solutions.
We refer to \cite{birge2011introduction} for an extensive overview on the research area of stochastic programming.

In summary, the STSP-GL combines the deterministic SGMLP and TSP-GL with stochastic demand as for an undirected graph with compulsory nodes, it aims to find a Hamiltonian tour on a subset of nodes, which includes all compulsory nodes, such that a given percentage of passenger demand is met with a given probability.

In recent years, researchers have made significant strides in integrating uncertain passenger demand into optimization models for public transportation systems. \cite{tian2021autonomous} consider the integration of autonomous vehicles into bus transit systems and determine the optimal bus fleet size and its assignment onto multiple bus lines in a bus service network considering uncertain demand. \cite{ma2021robust} study the problem of bus regulation with uncertain passenger demand and disturbance for urban rapid transit lines on the basis of a robust model predictive control algorithm. \cite{gkiotsalitis2019robust} incorporate travel time and passenger demand uncertainty to determine robust schedules that minimize the possible loss at worst-case scenarios. \cite{li2019robust} propose a robust dynamic control framework that considers delay disturbances and passenger demand uncertainty to reduce bus bunching. So far, uncertainty in passenger demand that allows a DAS operator to determine a general order of transit stops that is robust against volatile passenger demand has not been considered explicitly in the tactical planning phase of DAS. Accordingly, DAS operators are unable to analyze the impact of different service levels and the frequency of meeting them on the route design and served passenger demand.

\subsection{Contribution and Organization}
To close the research gap outlined above, we propose an algorithmic framework to determine a Hamiltonian tour on a subset of nodes of an undirected graph, which minimizes a linear combination of route design and passenger routing costs. The subset of nodes includes a set of compulsory nodes as well as other non-compulsory nodes, which have to be selected such that serving the passengers requiring transportation between nodes in this set meets a percentage of the uncertain passenger demand with a given probability. Accordingly, our problem consists of selecting a suitable subset of nodes and determining a Hamiltonian tour on it. Specifically, we formulate the so-called stochastic TSP with generalized latency (STSP-GL) as a chance-constrained program and determine its deterministic equivalent. To solve this deterministic equivalent efficiently, we present a branch-and-price (B\&P) approach, which allows us to decompose the STSP-GL into multiple TSP-GLs. Additionally, we develop a local search-based heuristic and a hybrid solution approach. The hybrid solution approach utilizes our heuristic approach in our B\&P approach. We show the efficiency of determining good bounds through our approaches by applying them to the TSPGL2 instance set and show that our three approaches, i.e., B\&P, heuristic, and hybrid, meet the MIP benchmark's bounds on small instances and find much better bounds on larger instances. Furthermore, we show that our three approaches find good upper bounds much faster than our MIP benchmark. Additionally, we show that when we compare our stochastic approach to a deterministic approach, as expected, the deterministic approach's solutions have lower design costs but do not meet the necessary feasibility frequencies, i.e., they are not feasible stochastic solutions. Finally, we show that the service levels influence the design costs much more than the feasibility frequencies. In summary, our work allows the integration of stochastic demand into the TSP-GL, which grants us more robust routes compared to the deterministic demand setting.

The remainder of this paper is as follows. We specify our problem setting in Section~\ref{sec: Problem Setting} and develop our methodology in Section~\ref{sec: Methodology}. In Section~\ref{sec: Results}, we describe a case study based on the TSPGL2 dataset, which we use to present numerical results that show the efficiency of our algorithmic framework. Section~\ref{sec: Conclusion} concludes this paper by summarizing its main findings.

\section{Problem description and mathematical formulations}
\label{sec: Problem Setting}
\begin{figure}[!t]
\centering
\begin{subfigure}[t]{0.45\textwidth}
  \centering
  \includegraphics[width=\columnwidth]{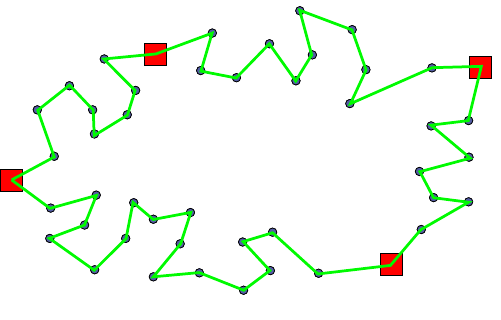}
    \captionsetup{format=hang}
	\caption{Solution for deterministic passenger demand (TSP-GL) with all optional stops}
\end{subfigure}
\begin{subfigure}[t]{0.45\textwidth}
  \centering
  \includegraphics[width=\columnwidth]{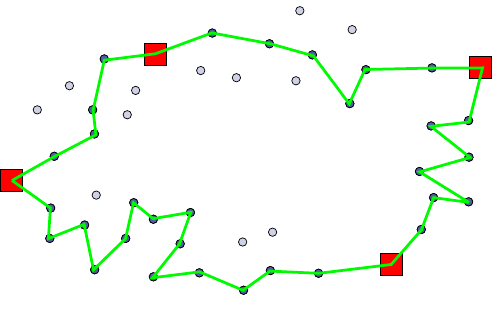}
    \captionsetup{format=hang}
    \caption{Solution for stochastic passenger demand (STSP-GL) with reduced sets of optional stops}
\end{subfigure}
\caption{Main components of the DAS lines for deterministic and stochastic passenger demand. Picture modified from \cite{crainic2012designing}}
\label{fig: Comparison TSP-GL and STSP-GL}
\end{figure}

In this work, we focus on activity (2) of a DAS line's aforementioned tactical planning activities, i.e., determining a general order of transit stops.
To do so, we build on the work of \cite{errico2017benders} by incorporating stochastic passenger demand into our problem setting, i.e., we assume the volume of demand and the presence of demand between two stops to be uncertain. Accordingly, we aim to choose a set of arcs to form the tour of a DAS before the passenger demand is known. We define a \textit{service level} as the percentage of passenger demand transported, which allows us to determine the order of compulsory stops and the service area, i.e., the sets of optional stops that must be covered between consecutive compulsory stops, such that they fulfill a given service level with a specified probability, e.g., 80\% of the estimated passenger demand has to be fulfilled with a probability of 90\% (see Figure~\ref{fig: Comparison TSP-GL and STSP-GL}). 

As a consequence, the STSP-GL involves three levels of decisions. First, we select a subset of transit stops from all possible transit stops such that the expected sum of passengers transported satisfies a specified service level of passenger demand with a predefined probability. This subset of transit stops must include all compulsory stops. Second, we choose a set of edges that form a Hamiltonian circuit on the selected subset of transit stops. Third, we determine the direction of passenger routes for the stochastic demand on the Hamiltonian circuit. In this tactical planning problem, no decisions are taken after uncertainty is resolved, and perfect information is known. Our objective is to minimize a combination of the design cost of the route, i.e., the sum of selected edge costs, as well as the expected average passenger routing cost.

\subsection{Chance-constrained STSP-GL}
We present a similar mathematical formulation to \cite{errico2017benders} to stay consistent with the literature. We consider a complete mixed graph $G = (N, E \cup A)$, where $N = \{1, \dots, n\}$ describes a set of nodes, $E = \{[i,j]: i,j \in N, i<j\}$ describes a set of edges that models an ordered tuple of nodes, and ${A = \{(i,j): i,j \in N, i \neq j\}}$ describes a set of undirected arcs. Here, edges model vehicle movement while arcs model passenger movement. We associate a \textit{design} cost $\overline{c}_{[i,j]} > 0$ to each edge $[i,j] \in E$. Moreover, let $d^{hk}(\xi) \geq 0$ be the uncertain amount of demand for each origin and destination node pair $(h,k)$. Here, $d^{hk}(\xi) \in \Xi$ describes a discrete random variable. Let $D = \{(h,k):h,k \in N, P(d^{hk}(\xi) > 0) > 0\}$ be the set of node pairs with a positive probability of non-zero demand. We associate a travel time $c_{ij}$ to each arc $(i,j) \in A$ and a random variable for the \textit{routing cost} $q_{ij}^{hk}(\xi) = c_{ij} \frac{d^{hk}(\xi)}{\theta \sum_{(s,t) \in D} d^{st}(\xi)}$ to every arc $(i,j) \in A$ and every node pair $(h,k) \in D$ that depends in turn on the random variable $d^{hk}(\xi)$. Furthermore, we set a \textit{service level} $\theta \in [0,1]$ and a \textit{feasibility frequency} $1 - \rho \in [0,1]$. Let $\alpha$ be a tradeoff parameter between design and routing costs. We seek to find a subset $T \subset E$ of edges, such that $T$ describes a tour in $G$ which visits all compulsory nodes $\mathcal{C}$ and enables $\theta$ percentage of all customers to travel from their origins to their destinations on the arcs whose corresponding edges are in the tour with a probability of $1 - \rho$.

Accordingly, we introduce decision variables $x_{[i,j]}$ which indicate whether edge $[i,j] \in E$ lies in $T$ ($x_{[i,j]} = 1$) or not ($x_{[i,j]} = 0$). Additionally, decision variables $f_{ij}^{hk}$ indicate the percentage of demand $(h,k) \in D$ on arc $(i,j) \in A$. Furthermore, decision variables $z^{hk}$ indicate whether tour $T$ serves demand pair $(h,k)$ ($z^{hk} = 1$) or not ($z^{hk} = 0$) and $w_i$ indicates whether stop $i \in N$ is visited by tour $T$ ($w_i = 1$) or not ($w_i = 0$). 
Consequently, $\sum_{(i,j) \in A} \sum_{(h,k) \in D} q_{ij}^{hk}(\xi)f_{ij}^{hk}$ measures the expected average passenger travel time for service level $\theta$.
Notice that, similar to the TSP-GL, we do not consider any capacity restrictions in our problem. We formulate the STSP-GL as follows:

{\small
\allowdisplaybreaks
    \setlength{\abovedisplayskip}{0pt}
    \setlength{\abovedisplayshortskip}{0pt}
    \setlength{\belowdisplayskip}{0pt}
    \setlength{\belowdisplayshortskip}{0pt}
    \begin{subequations}
        \label{CC-SGMLP basic}
        {\setlength{\belowdisplayskip}{0pt}
    \begin{align}
        \min_{x, f} \quad (1 - \alpha) \sum_{[i,j] \in E} \overline{c}_{[i,j]}x_{[i,j]} +  \alpha \mathbb{E} \left[ \sum_{(i,j) \in A} \sum_{(h,k) \in D} q_{ij}^{hk}(\xi)f_{ij}^{hk} \right]
        \label{obj: func}
    \end{align}}
    \begin{align}
        \mathcal{P}\Big(\sum_{(h,k) \in D} z^{hk}d^{hk}(\xi) \geq \theta \sum_{(h,k) \in D}d^{hk}(\xi)\Big) \geq 1-\rho
        \label{cons: selective bound}
    \end{align}
    \begin{align}
        \sum_{(i,j)\in A}f_{ij}^{hk} - \sum_{(j',i)\in A}f_{j'i}^{hk} &=\,
    \begin{cases}
      z^{hk}& \text{if } i = h,\\
      -z^{hk}& \text{if } i = k,\\
      0& i\neq h,k
    \end{cases}
    &\forall i \in N, \; \forall (h,k) \in D        
    \label{cons:flow cons f}\\
        w_i &= 1 &
        \quad \forall i \in \mathcal{C}
        \label{cons: compulsory stops}\\
        \sum_{[i,j] \in E} x_{[i,j]} &= 2 w_i &
        \quad \forall i \in N
        \label{cons: delta 1}\\
        w_h &\geq z^{hk}&
        \quad \forall (h,k) \in D
        \label{cons: w h}\\
        w_k &\geq z^{hk}&
        \quad \forall (h,k) \in D
        \label{cons: w k}\\
        f_{ij}^{hk} &\geq 0&
        \quad \forall (i,j) \in A, \; \forall (h,k) \in D\\
        x_{[i,j]} - f_{ij}^{hk} &\geq 0 &
        \quad \forall [i,j] \in E, \; \forall (h,k) \in D
        \label{cons: f ij}\\
        x_{[i,j]} - f_{ji}^{hk} &\geq 0&
        \quad \forall [i,j] \in E, \; \forall (h,k) \in D 
        \label{cons: f ji}\\
        x_{[i,j]} &\in \{0,1\}&
        \quad \forall [i,j] \in E
        \label{cons: int x}\\
        z^{hk} &\in \{0,1\} &
        \quad \forall (h,k) \in D
        \label{cons: int z}\\
        w_i &\in \{0,1\}&
        \quad \forall i \in N
        \label{cons: int w}
    \end{align}
    \end{subequations}
}
The objective function~\eqref{obj: func} accounts for the sum of design costs and the sum of expected passenger routing costs. Constraint~\eqref{cons: selective bound} imposes that with a probability of $1 - \rho$, a service level of $\theta$ is met. The multicommodity flow balance constraints~\eqref{cons:flow cons f} ensure that if demand pair $(h,k)$ is served ($z^{hk} = 1$), one unit of flow is sent from $h$ to $k$. Here, similar to the TSP-GL formulation, we model the routing of passenger demand by unit-flow variables and let $q_{ij}^{hk}$ account for the application-specific interpretation of these routing costs~\citep{errico2017benders}. Constraints~\eqref{cons: compulsory stops} ensure that every compulsory stop is included in the tour and Constraints~\eqref{cons: delta 1}~-~\eqref{cons: w k} impose that if demand pair $(h,k)$ is served ($z^{hk} = 1$), the selected edges incident to nodes $h$ and $k$ must be exactly two. Constraints~\eqref{cons: f ij} and~\eqref{cons: f ji} enforce that any commodity $(h,k)$ can only travel on arcs $(i,j)$ and $(j,i)$ if the corresponding edge $[i,j]$ is in tour $T$. Note that our formulation does not include subtour elimination constraints as we focus on problem instances in which probability distributions induce a single connected component, i.e., there do not exist disjunct subsets of nodes such that passengers mostly travel only between nodes of the same subset. This is the case if $D$ is dense enough and consequently does not induce more than one connected component. For a more thorough explanation on the absence of subtour elimination constraints, we refer to~\cite{errico2017benders}. 

\subsection{Deterministic equivalent of the STSP-GL}
The computational tractability of Problem~\ref{CC-SGMLP basic} depends on characteristics of the uncertainty considered, i.e., whether it is possible to decouple decisions from random variables to transform probabilistic into deterministic constraints. Our problem is a chance-constrained optimization problem as an operator aims for a tour that fulfills a service level for a share instead of all passenger requests. In this setting, we study the problem's sample counterpart, which becomes computationally more tractable as it bears only a finite number of constraints. We now present a deterministic equivalent of our stochastic program, which ensures that the service level is above $\theta$ to a predefined feasibility frequency of $1 - \rho$. We introduce binary variables $Y_{s}$, which indicate whether service level $\theta$ is fulfilled ($Y_{s} = 1$) or not ($Y_{s} = 0$) for a given scenario $s \in S$ in the scenario set of passengers with passenger demands $s^{hk}$ and reformulate Problem~\ref{CC-SGMLP basic} as Problem~\ref{CC-SGMLP decomp}.

{\small
\allowdisplaybreaks
    \setlength{\abovedisplayskip}{0pt}
    \setlength{\abovedisplayshortskip}{0pt}
    \setlength{\belowdisplayskip}{0pt}
    \setlength{\belowdisplayshortskip}{0pt}

    \begin{subequations}
    \label{CC-SGMLP decomp}
    \begin{align}
        \min_{x, f} \quad (1 - \alpha) \sum_{[i,j] \in E} \overline{c}_{[i,j]}x_{[i,j]} + \alpha \sum_{(i,j) \in A} \sum_{(h,k) \in D} \tilde{q}_{ij}^{hk}f_{ij}^{hk} 
        \label{obj: func decomp}
    \end{align}
    \begin{align}
        \sum_{(h,k) \in D} z^{hk} s^{hk} & \geq \theta  Y_{s} \sum_{(h,k) \in D} s^{hk}&
        \quad \forall s \in S \label{cons:selective bound 1 dec}\\
        \sum_{s \in S} Y_{s} &\geq (1 - \rho) \vert S \vert
        \label{cons:selective bound 2 dec}\\
        \sum_{(i,j)\in A}f_{ij}^{hk} - \sum_{(j',i)\in A}f_{j'i}^{hk} &=\,
    \begin{cases}
      z^{hk}& \text{if } i = h,\\
      -z^{hk}& \text{if } i = k,\\
      0& i\neq h,k
    \end{cases}
    &\forall i \in N, \; \forall (h,k) \in D
    \label{cons:flow cons f dec}\\
        w_i & = 1 &
        \quad \forall i \in \mathcal{C}
        \label{cons: compulsory stops dec}\\
        \sum_{[i,j] \in E} x_{[i,j]} & = 2 w_i &
        \quad \forall i \in N
        \label{cons: delta 1 dec}\\
        w_h & \geq z^{hk}&
        \quad \forall (h,k) \in D
        \label{cons: w h dec}\\
        w_k &\geq z^{hk}&
        \quad \forall (h,k) \in D
        \label{cons: w k dec}\\
        f_{ij}^{hk} &\geq 0&
        \quad \forall (i,j) \in A, \; \forall (h,k) \in D\\
        x_{[i,j]} - f_{ij}^{hk} & \geq 0 &
        \quad \forall [i,j] \in E, \; \forall (h,k) \in D
        \label{cons: f ij dec}\\
        x_{[i,j]} - f_{ji}^{hk} &\geq 0&
        \quad \forall [i,j] \in E, \; \forall (h,k) \in D 
        \label{cons: f ji dec}\\
        x_{[i,j]} &\in \{0,1\}&
        \quad \forall [i,j] \in E
        \label{cons: int x dec}\\
        z^{hk} &\in \{0,1\} &
        \quad \forall (h,k) \in D
        \label{cons: int z dec}\\
        w_i &\in \{0,1\}&
        \quad \forall i \in N
        \label{cons: int w dec}
    \end{align}
    \end{subequations}
    }

Here, we only reformulate Constraint~\eqref{cons: selective bound} as Constraints~\eqref{cons:selective bound 1 dec} and Constaint~\eqref{cons:selective bound 2 dec}. All other constraints are unchanged. Additionally, we update the passenger routing cost ${\tilde{q}_{ij}^{hk} = c_{ij} \Big(\frac{1}{\vert S\vert} \sum_{s \in S} \frac{s^{hk}}{\theta \sum_{(u,v) \in D} s^{uv}} \Big)}$ to account for the deterministic counterpart of our stochastic program. In his context, we assume the scenarios to be equally likely. 

\subsection{Decomposition of STSP-GL using feasibility covers}
This deterministic equivalent allows us to solve small instances but is still not tractable for larger ones since the number of combinations for $z$ is exponential in the size of $D$. Accordingly, we utilize the notion of feasibility covers to further decompose our model to utilize Benders decomposition \citep{errico2017benders} for the resulting subproblems.

\begin{definition}[Feasibility cover]\label{def: feasibility cover}
A \textit{feasibility cover} is a set $Q \subseteq D$ such that there exists a vector $\gamma \in \{0,1\}^{\vert S \vert}$ which is feasible for

\begin{subequations}
\begin{align}
 \sum_{(h,k) \in Q} s^{hk}   &\; \geq \; \theta  \gamma_{s} \sum_{(h,k) \in D} s^{hk} & \forall s \in S, \label{feasibility cover 1} \\
 \sum_{s \in S} \gamma_{s} &\; \geq \; (1 - \rho) \vert S \vert \label{feasibility cover 2}
\end{align}
\end{subequations}
\end{definition}

In other words, a feasibility cover is a subset of requests for which the service level is met for $\lceil (1 - \rho) \vert S \vert \rceil$ many scenarios.\smallskip

\begin{definition}[Minimal feasibility cover]\label{def: minimum feasibility cover}
A \textit{minimal feasibility cover} is a feasibility cover $Q \subseteq D$ such that there is no feasibility cover $Q’ \subseteq D$ such that $Q’ \subset Q$.
\end{definition}\smallskip

\begin{theorem}\label{theroem: minimal feasibility cover}
For an optimal solution $(x^*, f^*, z^*, Y^*)$ of Problem~\ref{CC-SGMLP decomp}, $K = \{(h,k): (z^{hk})^* = 1, (h,k) \in D\}$ is a minimal feasibility cover.
\end{theorem}

\proof
Assume we have an optimal solution $(x^*, f^*, z^*, Y^*)$ such that the corresponding set $K$ is not a minimal feasibility cover. 
Since $K$ is not a minimal feasibility cover, there has to exist a passenger request $(u,v) \in K$ for which there exists a $\gamma' \in \{0,1\}^{\vert S \vert}$  which fulfills 

\begin{subequations}
\begin{align}
 \sum_{(h,k) \in K-\{(u,v)\}} s^{hk}   &\; \geq \; \theta  \gamma_{s}' \sum_{(h,k) \in D} s^{hk} & \quad \forall \; s \in S\\
 \sum_{s \in S} \gamma_{s}' &\; \geq \; (1 - \rho) \vert S \vert &
\end{align}
\end{subequations}

Accordingly, ${(x', f', z', Y')}$ with 

\begin{align*}
&x' = x^*\\
&(f^{hk}_{ij})' = 
	\begin{cases}
      (f^{hk}_{ij})^*, & \text{if } (h,k) \neq (u,v),\\
      0, & \text{if } (h,k) = (u,v),\\
    \end{cases}\\
&(z^{hk})' = 
	\begin{cases}
      (z^{hk})^*, & \text{if } (h,k) \neq (u,v),\\
      0, & \text{if } (h,k) = (u,v),\\
    \end{cases}\\
&Y' = \gamma'
\end{align*}

is feasible for Problem \ref{CC-SGMLP decomp} and has an objective value of

\begin{align}
\small
    obj({(x', f', z', Y')}) = obj({(x^*, f^*, z^*, Y^*)}) - \nonumber\\
    \underbrace{\sum_{(i,j) \in A} c_{ij} \Big(\frac{1}{\vert S \vert} \sum_{s \in S} \frac{s^{hk}}{\theta \sum_{(r,t) \in D} s^{rt}} \Big) f_{ij}^{uv}}_{> 0}
\end{align}

This contradicts our assumption as ${(x^*, f^*, z^*, Y^*)}$ is an optimal solution. $\hfill \square$
\endproof
\smallskip

Definition~\ref{def: feasibility cover} allows us to reformulate Problem~\ref{CC-SGMLP decomp} by using feasibility covers. Furthermore, Definition~\ref{def: minimum feasibility cover} and Theorem~\ref{theroem: minimal feasibility cover} enable us to only consider minimal feasibility covers in the reformulation.

Let $\mathcal{X}_{\theta, \rho}$ be the set of all minimal feasibility covers of $D$. We can then formulate our original problem as Problem~\ref{cc-sgmlp reformulation}.
{\small
    \setlength{\abovedisplayskip}{0pt}
    \setlength{\abovedisplayshortskip}{0pt}
    \setlength{\belowdisplayskip}{0pt}
    \setlength{\belowdisplayshortskip}{0pt}
\begin{subequations}
\label{cc-sgmlp reformulation}

\begin{align}
    \min_{x, f} \quad (1 - \alpha) \sum_{[i,j] \in E} \overline{c}_{[i,j]}x_{[i,j]} + \alpha \sum_{(i,j) \in A} \sum_{(h,k) \in D} \tilde{q}_{ij}^{hk}f_{ij}^{hk} 
    \label{obj4: func}
\end{align}
\begin{align}
    \sum_{Q\in \mathcal{X}_{\theta, \rho}} \chi_Q & = 1
    \label{cons4: convex chi ref}\\
    \sum_{[i,j]\in E} x_{[ij]} - 2\sum_{Q\in \mathcal{X}_{\theta, \rho}} \chi_Q l_Q(i) & = 0 &
    \quad \forall i \in N
    \label{cons4: design max ref}\\
    \sum_{(i,j)\in A}f_{ij}^{hk} - \sum_{(j',i)\in A}f_{j'i}^{hk} &=
    \begin{cases}
      \sum_{Q\in \mathcal{X}_{\theta, \rho}} \chi_Q r_Q^{hk}, \; \text{if } i = h,\\
      -\sum_{Q\in \mathcal{X}_{\theta, \rho}} \chi_Q r_Q^{hk}, \;  \text{if } i = k,\\
      0, \text{else}
    \end{cases}
    &\quad\forall i \in N, \; \forall \; (h,k) \in D
    \label{cons4: flow cons f ref}\\
    f_{ij}^{hk} &\geq 0&
    \quad \forall (i,j) \in A, \; \forall (h,k) \in D 
    \label{cons4: passenger flow bound ref}\\
    x_{[i,j]} - f_{ij}^{hk} &\geq 0 &
    \quad \forall [i,j] \in E, \; \forall (h,k) \in D 
    \label{cons4: f ij}\\
    x_{[i,j]} - f_{ji}^{hk} &\geq 0&
    \quad \forall [i,j] \in E, \; \forall (h,k) \in D 
    \label{cons4: f ji}\\
    x_{ij} &\in \{0,1\} &
    \quad \forall [i,j] \in E
    \label{cons4: int x ref}\\
    \chi_{Q} &\in \{0,1\} &
    \quad \forall Q \in \mathcal{X}_{\theta, \rho} 
    \label{cons4: int chi ref}
\end{align}
\end{subequations}
}

Here, $l_Q(i) = 1$ if $i$ is a compulsory stop ($i \in \mathcal{C}$) or there exists a demand pair $(h,k) \in Q$ with $i = h$ or $i = k$ and $l_Q(i) = 0$ otherwise. Furthermore, $r_Q^{hk} = 1$ if ${(h,k) \in Q}$ and $r_Q^{hk} = 0$ otherwise. Constraint~\eqref{cons4: convex chi ref} ensures that exactly one minimal feasibility cover is chosen. Constraints~\eqref{cons4: design max ref} correspond to Constraints~\eqref{cons: delta 1 dec} and ensure that a tour visits all stops for which $l_Q(i) = 1$ if $Q$ is chosen in~\eqref{cons4: design max ref}. Constraints~\eqref{cons4: flow cons f ref} are flow conservation constraints for all requests in $Q$ if $Q$ is chosen in~\eqref{cons4: design max ref}. Constraints~\eqref{cons4: passenger flow bound ref}-\eqref{cons4: int x ref} did not change compared to Problem~\ref{CC-SGMLP basic}. Constraints~\eqref{cons4: int chi ref} define $\chi_{Q}$'s domain.

In the remainder of this section, we show that we can find an optimal solution of our STSP-GL by determining a minimal feasibilty cover whose corresponding TSP-GL (see Problem~\ref{TSP-GL}) has the smallest objective function out of all minimal feasibility covers' corresponding TSP-GLs. To do so, we first show the corresponding TSP-GL to a minimal feasibility cover (see Problem~\ref{TSP-GL}) can be solved on a subgraph $G'$ of our original graph $G$. Second, we show that once we chose a minimal feasibility cover in our STSP-GL, the STSP-GL formulation decomposes into the aformentioned TSP-GL (see Theorem~\ref{theorem: minimal cover optimal solution}).

Let us revisit the TSP-GL given a feasibility cover $Q$ on the complete mixed graph $G'(Q) \subset G$ where $G'(Q) = (N'(Q), E'(Q) \cup A'(Q))$ with $N'(Q) =  \{i : (i,j) \in Q \; or \; (j',i) \in Q\}$, $E'(Q) = \{[i,j] : [i,j] \in E, \; i,j \in N'(Q)\}$ and $A'(Q) = \{(i,j) : (i,j) \in A, \; i,j \in N'(Q)\}$. We formulate the corresponding TSP-GL for $Q$ as follows:
\smallskip
{\small
    \setlength{\abovedisplayskip}{0pt}
    \setlength{\abovedisplayshortskip}{0pt}
    \setlength{\belowdisplayskip}{0pt}
    \setlength{\belowdisplayshortskip}{0pt}
\begin{subequations}
\label{TSP-GL}

    \begin{align}
        \min_{\bar{x}, \bar{f}} (1 - \alpha) \sum_{[i,j] \in E'(Q)} \overline{c}_{[i,j]} \bar{x}_{[i,j]} + \alpha \sum_{(i,j) \in A'(Q)} \sum_{(h,k) \in Q} \tilde{q}_{ij}^{hk} \bar{f}_{ij}^{hk} 
        \label{obj: func TSP-GL}
    \end{align}
    \begin{align}
        \sum_{(i,j)\in A'(Q)} \bar{f}_{ij}^{hk} - \sum_{(j',i)\in A'(Q)} \bar{f}_{j'i}^{hk} &=\,
    \begin{cases}
      1& \text{if } i = h,\\
      -1& \text{if } i = k,\\
      0& i\neq h,k
    \end{cases}
    &\forall i \in N'(Q), \; \forall (h,k) \in Q
    \label{cons: flow cons f TSP-GL}\\
        \sum_{[i,j] \in E'(Q)} \bar{x}_{[i,j]} & = 2&
        \quad \forall i \in N'(Q)
        \label{cons: design TSP-GL}\\
        \bar{f}_{ij}^{hk} &\geq 0&
        \quad \forall (i,j) \in A'(Q), \; \forall (h,k) \in Q 
        \label{cons: passenger flow bound TSP-GL}
\end{align}
\begin{align}
        \bar{x}_{[i,j]} - \bar{f}_{ij}^{hk} &\geq 0&
        \quad \forall [i,j] \in E'(Q), \; \forall (h,k) \in Q 
        \label{cons: f ij TSP-GL}\\
        \bar{x}_{[i,j]} - \bar{f}_{ji}^{hk} &\geq 0&
        \quad \forall [i,j] \in E'(Q), \; \forall (h,k) \in Q 
        \label{cons: f ji TSP-GL}\\
        \bar{x}_{[i,j]} &\in \{0,1\}&
        \quad \forall [i,j] \in E'(Q)
        \label{cons: int x TSP-GL}
    \end{align}
\end{subequations}
}

\begin{theorem}\label{theorem: x and f zero}
If Problem~\ref{cc-sgmlp reformulation} is feasible and triangle inequality holds for $\bar{c}$, there exists an optimal solution, in which for each served request in an optimal minimal feasibility cover $Q^*$, i.e., $r \in Q^*$ with $\chi_{Q^*} = 1$, it holds that $\bar{x}_{[i,j]}^* = 0 \; \forall [i,j] \in E \setminus E'(Q^*)$ and $(\bar{f}_{ij}^{hk})^* = 0 \; \forall (i,j) \in A,\; \forall (h,k) \in D$ if $(i,j) \notin A'(Q^*)$ or $(h,k) \notin Q^*$.
\end{theorem}

\proof
All edges in $E \setminus E'(Q^*)$ are adjacent to a node $i \notin N'(Q^*)$, i.e., a node to which no request in $Q^*$ has an origin or destination. If the tour formed by $\bar{x}$ includes such an edge, the design costs of an alternative tour which skips all nodes $i \notin N'(Q^*)$ would always be lower or equal to the original route because triangle inequality holds. Similarly, as no request in $Q^*$ has an origin or destination in a node $i \notin N'(Q^*)$, alternative passenger routes could also skip these nodes, which, because of triangle inequality, leads to passenger routes with equal or smaller costs. As the original passenger routes only use arcs whose corresponding edges are in the tour formed by $\bar{x}$ and the alternative passenger routes would skip the same nodes as the alternative tour, all passenger routes would still be feasible. $\hfill \square$
\endproof\smallskip

\begin{definition}[Characteristic function of TSP-GL variables]\label{def: characteristic func TSP-GL}
For a solution $(\bar{x}, \bar{f})$ of Problem~\ref{TSP-GL} with minimal feasibility cover $Q \in \mathcal{X}_{\theta, \rho}$, define $\hat{x}(Q)$ and $\hat{f}(Q)$ as follows:
\begin{equation}
\hat{x}_{[i,j]}(Q) =
	\begin{cases}
      \bar{x}_{[i,j]}& \text{if } [i,j] \in E'(Q),\\
      0& \text{else}
    \end{cases}
    \quad
    \forall [i,j] \in E
\end{equation}

{\small
    \setlength{\abovedisplayskip}{0pt}
    \setlength{\abovedisplayshortskip}{0pt}
    \setlength{\belowdisplayskip}{-8pt}
    \setlength{\belowdisplayshortskip}{0pt}
\begin{align}
    \hat{f}_{i,j}^{h,k}(Q) =
	\begin{cases}
      \bar{f}_{i,j}^{h,k}& \text{if } (i,j) \in A'(Q) \; and \; (h,k) \in Q,\\
      0& \text{else}
    \end{cases}
    &\quad \forall (i,j) \in E \; and \; (h,k) \in D
\end{align}
}
\end{definition}

\vspace{2ex}
With the new STSP-GL formulation in Problem~\ref{cc-sgmlp reformulation}, the TSP-GL formulation for a minimal feasibility cover $Q$ in Problem~\ref{TSP-GL}, and the definition of $\hat{x}(Q)$ and $\hat{f}(Q)$, we can show the following theorem.\medskip

\begin{theorem}\label{theorem: minimal cover optimal solution}
If triangle inequality holds for $\bar{c}$ then for all $\chi = \mathbbm{1}_{\{Q = Q'\}}$, Problem~\ref{cc-sgmlp reformulation} decomposes into Problem~\ref{TSP-GL}. Furthermore, $(\mathbbm{1}_{\{Q = Q^*\}}, \hat{x}(Q^*), \hat{f}(Q^*))$, with $Q^* = arg\,min \{obj_{TSP-GL}(Q) \,:\, Q \in \mathcal{X}_{\theta, \rho}\}$ is an optimal solution of Problem~\ref{cc-sgmlp reformulation}.
\end{theorem}\smallskip

\proof
According to Therorem~\ref{theorem: x and f zero}, we can assume that for an optimal solution $(\mathbbm{1}_{\{Q = Q^{**}\}}, \bar{x}(Q^{**}), \bar{f}(Q^{**}))$ it holds that $\bar{x}_{[i,j]}^{**} = 0 \; \forall [i,j] \in E \setminus E'(Q^{**})$ and $(\bar{f}_{ij}^{hk})^{**} = 0 \; \forall (i,j) \in A,\; \forall (h,k) \in D$ if $(i,j) \notin A'(Q^{**})$ or $(h,k) \notin Q^{**}$. Thus, we can fix these decision variables for each $Q' \in \mathcal{X}_{\theta, \rho}$. In this case, the objective values of~\eqref{obj4: func} and~\eqref{obj: func TSP-GL} are equal. 
Furthermore, since $\chi = \mathbbm{1}_{\{Q = Q'\}}$, Constraints~\eqref{cons4: design max ref} are equal to Constraints~\eqref{cons: design TSP-GL} $\forall i \in N'(Q')$.
By definition of $r_{Q'}^{hk}$ and since $\chi_{Q'} = 1$, Constraints~\eqref{cons4: flow cons f ref} are equal to Constraints~\eqref{cons: flow cons f TSP-GL} $\forall i \in N'(Q'), \; \forall (h,k) \in Q'$.
Additionally, Constraints~\eqref{cons4: passenger flow bound ref} - \eqref{cons4: int x ref} are equal to Constraints~\eqref{cons: passenger flow bound TSP-GL}~-~\eqref{cons: int x TSP-GL} for all variables which are not fixed through Theorem~\ref{theorem: x and f zero}.
Lastly, for each $Q' \in \mathcal{X}_{\theta, \rho}$, Constraint~\eqref{cons4: convex chi ref} and Constraints~\eqref{cons4: int chi ref} hold for $\chi = \mathbbm{1}_{\{Q = Q'\}}$.
This shows, that if the triangle inequality holds for $\bar{c}$ then for all $\chi = \mathbbm{1}_{\{Q = Q'\}}$, Problem~\ref{cc-sgmlp reformulation} decomposes into Problem~\ref{TSP-GL}.

The second part of the theorem follows directly. $\hfill \square$
\endproof
\section{Methodology}
\label{sec: Methodology}
In this section, we first show a standard B\&P approach (Section~\ref{subsec: BP}) on which we base our adapted B\&P approach (Section~\ref{subsec: Adapted BP}). Then, we introduce a heuristic approach that allows us to solve the STSP-GL without finding solutions for its continuous relaxation (Section~\ref{subsec: heuristic}). Finally, we propose a hybrid approach in which we combine our adapted B\&P and heuristic approaches (Section~\ref{subsec: hybrid}).

\subsection{Branch-and-price approach}
\label{subsec: BP}
\begin{figure}[!t]
  \centering
  \includegraphics[width=0.4\textwidth]{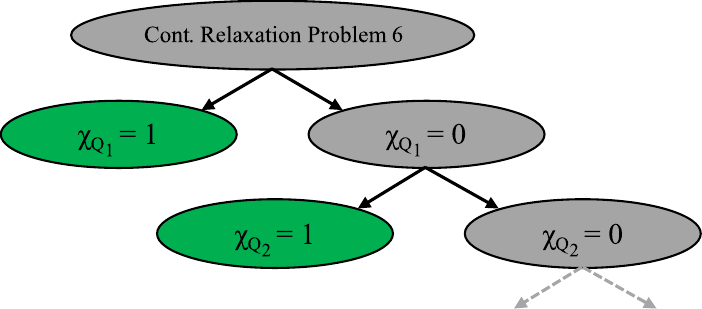}
	\caption{Branch-and-Price algorithm}
	\label{fig:branch-and-price}
\end{figure}

According to Theorem~\ref{theorem: minimal cover optimal solution}, we must find the minimal feasibility cover $Q^*$ to solve Problem~\ref{cc-sgmlp reformulation}. Even determining $\mathcal{X}_{\theta, \rho}$ is not a trivial task as the number of minimal feasibility covers may grow exponentially in the size of all requests $D$. One straightforward approach to find $Q^*$ is to embed Problem~\ref{cc-sgmlp reformulation} in a B\&P framework (Figure~\ref{fig:branch-and-price}) in which the nodes in the branch-and-bound tree (green) with integer decision variables $x_{[i,j]}$ decompose into TSP-GLs.

In such a B\&P framework, we start with a small subset $\mathcal{X}_{\theta, \rho}'$ of minimal feasibility covers and iteratively add new minimal feasibility covers to $\mathcal{X}_{\theta, \rho}'$ via column generation by solving the so-called pricing problem.
In this context, let $\iota_i \in \mathbbm{R}$ and $\epsilon_h^{hk}, \epsilon_k^{hk} \in \mathbbm{R}$ be the dual variables of Constraints~\eqref{cons4: design max ref} and \eqref{cons4: flow cons f ref} respectively. The scalar $\beta$ is the dual value of Constraint~\eqref{cons4: convex chi ref}. We aim to find a minimal feasibility cover $Q$ in each iteration of our column generation for which the reduced cost condition \eqref{reduced cost formula}, that indicates whether a variable should be brought into the basis, holds.

\begin{gather}
2 \sum_{i \in N} \iota_i l_Q(i) + \sum_{(h,k) \in D} r_Q^{hk} (\epsilon_{h}^{hk} - \epsilon_{k}^{hk}) - \beta < 0
\label{reduced cost formula}
\end{gather}

Accordingly, we formulate our pricing problem as follows.
{\small
    \setlength{\abovedisplayskip}{1pt}
    \setlength{\abovedisplayshortskip}{1pt}
    \setlength{\belowdisplayskip}{1pt}
    \setlength{\belowdisplayshortskip}{1pt}
\begin{subequations}
\label{cc-sgmlp pricing ref}
    \begin{align}
        \min_{l, r} \quad 2\sum_{i \in N} \iota_i l_i + \sum_{(h,k) \in D} r^{hk} (\epsilon^{hk}_h - \epsilon^{hk}_k) - \beta 
        \label{obj: func pricing ref}
    \end{align}
    \begin{align}
        l_i & = 1 &
        \quad \forall \; i \in \mathcal{C} 
        \label{cons: pricing compulsory stops}\\
        l_h & \geq r^{hk} &
        \quad \forall \; (h,k) \in D 
        \label{cons: pricing ref l_h}\\
        l_k & \geq r^{hk} &
        \quad \forall \; (h,k) \in D  
        \label{cons: pricing ref l_k}\\
        \sum_{(h,k) \in D} r^{hk} s^{hk} & \geq \gamma_{s} \theta \sum_{(h,k) \in D} s^{hk} &
        \quad \forall \; s \in S  
        \label{cons: pricing ref service level}\\
        \sum_{s\in S} \gamma_{s} & \geq (1 - \rho) \vert S \vert
        \label{cons: pricing ref feasibility target}\\
        l_i &\in \{0,1\} &
        \quad \forall \; i \in N 
        \label{cons: pricing ref l int}\\
        r^{hk} &\in \{0,1\} &
        \quad \forall \; (h,k) \in D 
        \label{cons: pricing ref q int}\\
        \gamma_{s} &\in \{0,1\} &
        \quad \forall \; s \in S 
        \label{cons: pricing ref gamma int}
    \end{align}
\end{subequations}
}

In Problem~\ref{cc-sgmlp pricing ref}, the decision variables $r^{hk}$ determine whether request $(h, k) \in D$ is in the feasibility cover ($r^{hk} = 1$) or not ($r^{hk} = 0$). The decision variables $l_i$ indicate whether at least one request $(i, k)$ or $(h, i)$ is in the feasibility cover or $i$ is a compulsory node ($l_i = 1$) or not ($l_i = 0$). This pricing problem does not ensure that the feasibility covers determined by $r^{hk}$ are minimal. We can enforce this by adding Constraints~\eqref{minimal feasibility covers constraint} as lazy constraints. Here, $\Psi$ is the set of feasibility covers that are not minimal.
\begin{gather}
\sum_{(h,k) \in Q} r^{hk}\leq \vert Q \vert -1 \qquad \forall Q \in \Psi
\label{minimal feasibility covers constraint}
\end{gather}

To consider the branching constraints $\chi_{Q_i} = 0$ of the B\&P, we extend the pricing problem by Constraints~\eqref{feasibility cover cut}.
\begin{gather}\label{feasibility cover cut}
\sum_{(h,k) \in Q_i} r^{hk} \leq \vert Q_i \vert -1 \qquad \forall Q_i \in \Phi
\end{gather}

Here, $\Phi$ is the set of feasibility covers $Q_i$ on which we already branched, i.e., for which constraints $\chi_{Q_i} = 0$ exists. 

\subsection{Adapted branch-and-price approach}
\label{subsec: Adapted BP}
One drawback of such a straightforward approach is that the large number of minimal feasibility covers makes solving even the root node of the B\&P approach very hard. Accordingly, we propose an adapted B\&P approach in which we iteratively find a set of feasibility covers through our Pricing Problem~\ref{cc-sgmlp pricing ref} without solving the continuous relaxation of our problem to optimality and solve their corresponding TSP-GLs to determine feasible solutions.

\begin{figure}[!t]
  \centering
  \includegraphics[width=\textwidth]{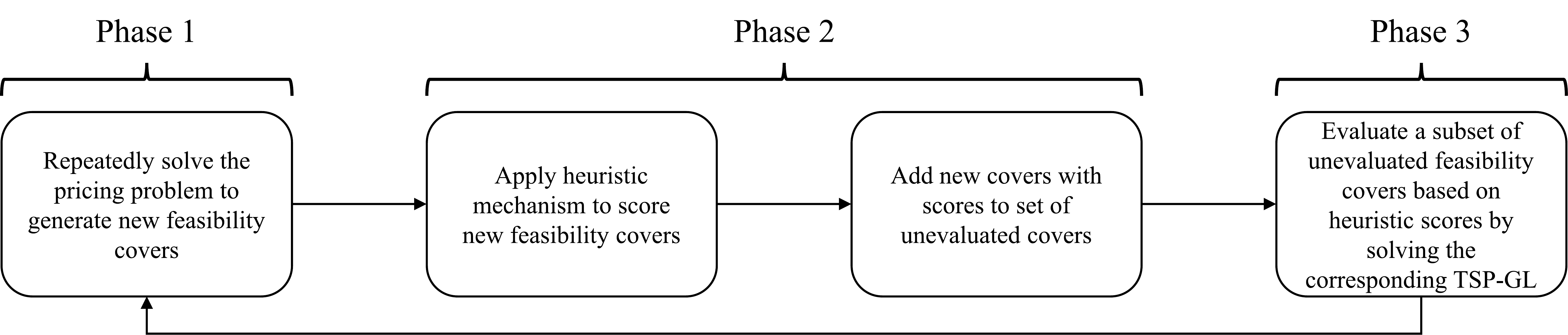}
	\caption{Flowchart B\&P approach algorithm}
	\label{fig:Flowchart SSTSP-GL}
\end{figure}

We initialize our algorithm by building our restricted master problem (RMP) with the empty set of minimal feasibility covers and variable $\chi_D$, i.e., the feasibility cover, which includes all requests. Additionally, we build our Pricing Problem~\ref{cc-sgmlp pricing ref} with an empty set of branching Constraints~\eqref{feasibility cover cut}. Afterward, we divide each iteration of our algorithm into three phases (see Figure~\ref{fig:Flowchart SSTSP-GL}). 

In the first phase, we repeatedly solve Pricing Problem~\ref{cc-sgmlp pricing ref} as an IP to generate new minimal feasibility covers. Whenever we find a new minimal feasibility cover with negative reduced cost, we update the lower bound of our B\&P and add it to the RMP. 

In the second phase, we apply a heuristic mechanism to score these new found minimal feasibility covers and add them to a list of unevaluated covers. By doing so, we allow for an efficient evaluation of feasibility covers, as solving the corresponding TSP-GL of a minimal feasibility cover takes relatively long, and our heuristic mechanism allows us to discard bad minimal feasibility covers much faster than solving the corresponding TSP-GL. Our heuristic mechanism scores a minimal feasibility cover $Q$ by evaluating its upper and lower bounds for the TSP-GL.  
If, for a minimal feasibility cover $Q$, the lower bound of its corresponding TSP-GL is higher than the upper bound of our B\&P, we can discard $Q$ as an optimal solution.

In the third phase, we evaluate a subset of unevaluated minimal feasibility covers chosen based on the heuristic scores by solving their corresponding TSP-GLs. Accordingly, we branch on the most promising minimal feasibility covers, i.e., the minimal feasibility covers in the set of unevaluated covers with the smallest upper bounds.

\subsubsection*{Determination of upper and lower bounds:}
In the second phase of our B\&P approach, we score feasibility covers. Here, the objective of the TSP-GL contains two potentially conflicting terms, i.e., design cost $\sum_{[i,j] \in E} \overline{c}_{[i,j]}x_{[i,j]}$ and routing cost $\sum_{(i,j) \in A} \sum_{(h,k) \in D} \tilde{q}_{ij}^{hk}f_{ij}^{hk}$. To derive a lower bound on the objective, we derive lower bounds on the individual terms. Accordingly, we solve the symmetric TSP on compulsory nodes and nodes with requests going from or to them in $Q$ to obtain a lower bound on the design costs ($lb_{design}$). Furthermore, for each request $(h, k) \in Q$, we assume a lower bound of $\tilde{q}_{hk}^{hk}$ for the routing cost ($lb_{routing}^{hk}$), which corresponds to a direct route from the origin $h$ to the destination $k$. Accordingly, $(1 - \alpha) lb_{design} + \alpha \sum_{(h, k) \in Q} lb_{routing}^{hk}$ is a lower bound for the TSP-GL on~$Q$. For the upper bound, we compute the shortest path of length $ub_{routing}^{hk}$ for each request $(h, k) \in Q$ on the solution of the symmetric TSP. Accordingly, $(1 - \alpha) lb_{design} + \alpha \sum_{(h, k) \in Q} ub_{routing}^{hk}$ is an upper bound for the TSP-GL on $Q$.

\subsubsection*{Benders decomposition for the TSP-GL:}
To solve the corresponding TSP-GLs in the third phase of our B\&P approach, similarly to \cite{errico2017benders}, we observe that for fixed design decision variables $\bar{x}$, the TSP-GL reduces to solving a one-commodity minimum cost flow model for each request in our minimal feasibility cover. Accordingly, we utilize a Benders decomposition approach to solve the TSP-GL. Denoting the routing cost of request $(h, k) \in Q$ with fixed design decision variables $\bar{x}$ as $z^{hk} (\bar{x})$, we define our PrimalSubproblem (PS) for each request $(h, k)$ as   

{\small
    \setlength{\abovedisplayskip}{0pt}
    \setlength{\abovedisplayshortskip}{0pt}
    \setlength{\belowdisplayskip}{0pt}
    \setlength{\belowdisplayshortskip}{0pt}
\begin{subequations}
\label{PS}
\begin{align}
    z^{hk}(\bar{x}) = \min_{\bar{f}} \quad \sum_{(i,j) \in A'(Q)} \tilde{q}_{ij}^{hk} \bar{f}_{ij}^{hk} \label{obj: func PS}
\end{align}
\begin{align}
    \sum_{(i,j)\in A'(Q)} \bar{f}_{ij}^{hk} - \sum_{(j',i)\in A'(Q)} \bar{f}_{j'i}^{hk} &=\,
    \begin{cases}
      1& \text{if } i = h,\\
      -1& \text{if } i = k,\\
      0& i\neq h,k
    \end{cases}
    &\forall i \in N'(Q)    
    \label{cons: flow cons f PS}\\
    \bar{f}_{ij}^{hk} & \geq 0 &
    \quad \forall (i,j) \in A'(Q), \; \forall (h,k) \in Q 
    \label{cons: passenger flow bound PS}\\
    \bar{x}_{[i,j]} - \bar{f}_{ij}^{hk} & \geq 0 & 
    \quad \forall [i,j] \in E'(Q), \; \forall (h,k) \in Q 
    \label{cons: f ij PS}\\
    \bar{x}_{[i,j]} - \bar{f}_{ji}^{hk} & \geq 0 &
    \quad \forall [i,j] \in E'(Q), \; \forall (h,k) \in Q 
    \label{cons: f ji PS}
\end{align}
\end{subequations}
}

As strong duality holds for each flow subproblem, we consider the DualSubproblems (DS) for every flow subproblem.

{\small
    \setlength{\abovedisplayskip}{0pt}
    \setlength{\abovedisplayshortskip}{0pt}
    \setlength{\belowdisplayskip}{0pt}
    \setlength{\belowdisplayshortskip}{0pt}
    \begin{subequations}
\label{DS}
\begin{align}
z^{hk}(\bar{x}) = \max \quad p_h^{hk} - p_k^{hk} - \sum_{[i,j] \in E'(Q)}\bar{\lambda}_{[i,j]}^{hk} \bar{x}_{[i,j]} 
\label{obj: DS}
\end{align}
\begin{align}
    p_i^{hk} - p_j^{hk} -\lambda_{ij}^{hk} &
    \leq q_{ij}^{hk} &
    \quad \forall (i,j) \in A'(Q) 
    \label{cons DS}\\
    \lambda_{ij}^{hk} & \geq 0 &
    \quad \forall (i,j) \in A'(Q) 
    \label{cons: lambda DS}\\
    p_{i}^{hk} & \geq 0 &
    \quad \forall (i, j) \in A'(Q) 
    \label{cons: p DS}
\end{align}
\end{subequations}
}
 
 Here, we define $\bar{\lambda}_{[i,j]} := \lambda_{ij} + \lambda_{ji}, \; \forall [i, j] \in E'(Q)$. We now denote $extr(R)$ and $rays(R)$ as the set of extreme points and rays of $R = \{(\lambda, p) \in \mathbbm{R}^{(\vert A'(Q) \vert + \vert N'(Q) \vert) \vert Q \vert} \vert (\lambda, p) \text{ that satisfy } \eqref{cons DS} - \eqref{cons: p DS}\}$ respectively.
 Consequently, we can express Problem~\ref{TSP-GL} as Problem~\ref{TSP-GL ref}.

{\small
    \setlength{\abovedisplayskip}{0pt}
    \setlength{\abovedisplayshortskip}{0pt}
    \setlength{\belowdisplayskip}{0pt}
    \setlength{\belowdisplayshortskip}{0pt}
    \begin{subequations}
\label{TSP-GL ref}
\begin{align}
\min \quad \alpha \eta + (1 - \alpha) \hspace{-10pt} \sum_{[i,j] \in E'(Q)} \bar{c}_{[i,j]} x_{[i,j]} 
\label{obj: TSP-GL ref}
\end{align}
\begin{align}
    \eta + \hspace{-10pt}\sum_{\substack{[i,j] \in E'(Q)\\
    (h, k) \in Q}}
    \bar{\lambda}_{[i,j]}^{hk} x_{[i,j]} & \geq \sum_{(h, k) \in Q} (p_h^{hk} - p_k^{hk}) &
    \forall (\lambda, p) \in extr(R) 
    \label{cons: optimality TSP-GL ref}\\
    \sum_{\substack{[i,j] \in E'(Q)\\
    (h, k) \in Q}}
    \bar{\lambda}_{[i,j]}^{hk} x_{[i,j]} & \geq \sum_{(h, k) \in Q} (p_h^{hk} - p_k^{hk}) & 
    \forall (\lambda, p) \in rays(R) 
    \label{cons: feasibility TSP-GL ref}\\
    \sum_{[i,j] \in E'(Q)} x_{[i,j]} & = 2&
    \forall i \in N'(Q)
    \label{cons: delta ref tsp-gl}\\
    x_{[i,j]} & \in \{0,1\} & 
    \forall [i,j] \in E'(Q)
    \label{cons: int x ref tsp-gl}
\end{align}
\end{subequations}
}

To solve Problem~\ref{TSP-GL ref}, we have to determine a set of feasibility cuts \eqref{cons: feasibility TSP-GL ref} and optimality cuts \eqref{cons: optimality TSP-GL ref}. We initialize our algorithm by using the solution of the symmetric TSP, which we solved for the bounds of $Q$, which is equivalent to initialization \textbf{h2} in \cite{errico2017benders}. We utilize a subtour elimination constraint strategy for our feasibility cut generation. As such, we adopt \eqref{feasibility cut benders} as feasibility cuts. Here $\Theta$ is the set of all subtours in the current solution.
\begin{align}
\label{feasibility cut benders}
\footnotesize
\sum_{[i,j] \in S} x_{[i,j]} \leq \vert S \vert -1 \qquad \forall S \in \Theta
\end{align}

According to \cite{errico2017benders}, we can add cuts \eqref{aggregated optimality cuts} as optimality cuts.
\begin{align}
 \footnotesize
\hspace{-5pt} \eta^{h} + \sum_{\substack{k \in N \textbackslash [h]\\
e \in E}}
\bar{\lambda}_e^{hk} x_e \geq \sum_{k \in N \textbackslash [h]} p_h^{hk} - p_k^{hk} \qquad \forall h \in N'(Q)
\label{aggregated optimality cuts}
\end{align}

Here, $\eta^h$ is the contribution, of all commodities with origin $h$, with $\sum_{h \in N} \eta^{h} = \eta$, to $\eta$.

Accordingly, in our Benders approach, we initialize Problem~\ref{TSP-GL ref} with a set of feasibility cuts, which are equal to the subtour elimination constraints we generated in the calculation of our bounds for the minimal feasibility covers and an empty set of optimality cuts. We additionally use the solution of the symmetric TSP obtained in the calculation of the bounds to generate optimality cuts. We then resolve Problem~\ref{TSP-GL ref} and check if we can add new feasibility cuts. If not, we solve Problem~\ref{DS} for each request $(h,k) \in Q$ and add new optimality cuts. We begin a new iteration by resolving Problem~\ref{TSP-GL ref} and terminate when we can not find new feasibility or optimality cuts. Note that we can add a set of feasibility cuts based on all found subtours in the calculation of the symmetric TSP when we initialize our Benders approach.

\subsection{Local search-based heuristic approach}
\label{subsec: heuristic}
One primary concern with this B\&P approach is the time required to solve the RMP in our column generation as the number of flow variables $f_{ij}^{hk}$ and the number of Constraints~\eqref{cons4: flow cons f ref}, \eqref{cons4: f ij}, and \eqref{cons4: f ji} grow by a factor of $\mathcal{O}(n^4)$ in the size of the underlying graph. Accordingly, we develop a local search-based heuristic that determines minimal feasibility covers without having to solve the RMP by leveraging the notion of a minimal node cover.

\begin{definition}[Minimal node cover]\label{def: minimum node cover}
We define a \textit{minimal node cover} $V'$ to be a set for which it holds that $\mathcal{C} \subseteq V' \subseteq V$ and there exists a feasibility cover in $\{(h,k) \, : \, h, k \in V'\}$ and no feasibility cover in ${\{(h,k): h, k \in V'\setminus \{v\}\}}$ for all $v \in V'$.
\end{definition}

Our heuristic consists of four main parts. An exploration step, a local search, a heuristic mechanism to score minimal feasibility covers, and the solving step. In the exploration step, the heuristic aims to find new minimal feasibility covers based on minimal node covers. To do so, we start with a random minimal node cover $V' \subseteq V$ with a predetermined size $\vert V' \vert = n$, which has not been found so far. We then determine a minimal feasibility cover $Q \subseteq Q' := \{(h,k) \, : \, h, k \in V'\}$. Algorithm~\ref{alg: get_minimum_covers} outlines how to find a minimal feasibility cover from a feasibility cover $Q'$. Note that if no feasibility cover of size $|Q| - 1$ exists for a feasibility cover $Q$, $Q$ is minimal. Utilizing this, we initialize Algorithm~\ref{alg: get_minimum_covers} with a feasibility cover $Q'$ with the aim of determining a minimal feasibility cover $Q \subseteq Q'$. To do so, we assign $Q \gets Q'$ (l. 1) and check whether there exists a feasibility cover $\bar{Q} \subset Q$ with $|\bar{Q}| = |Q| - 1$ (l.2 - 6). If that is the case, we recursively call Algorithm~\ref{alg: get_minimum_covers} with $\bar{Q}$ as an input (l.5). Otherwise, we return the minimal feasibility cover $Q$. 

\begin{algorithm}[!t]
\caption{Determination of minimal feasibility cover}
\label{alg: get_minimum_covers}
\KwData{feasibility cover $Q'$}
\KwResult{minimal feasibility cover $Q$}
$Q \gets Q'$\\
  \For{$r \in Q$}{
    \If{$Q \setminus \{r\} \text{feasibility cover}$}{
        $\bar{Q} \gets Q \setminus \{r\}$\\
        \Return minimal\textunderscore feasibility\textunderscore cover($\bar{Q}$)
    }
  }
  \Return $Q$
\end{algorithm}

\begin{algorithm}[!t]
\caption{Local search}
\label{alg: local search}
\KwData{minimal feasibility cover $Q'$; all requests $D$}
\KwResult{minimal feasibility cover $Q$}

$Q \gets Q'$\\
$R \gets \emptyset$

\While{$Q$ is feasibility cover}{
$R.add(Q.pop())$
}

\For{$r \in D \setminus R$}{
    $Q.add(r)$\\
    \If{$Q$ is feasibility cover}{
    break
    }
}
\If{$Q$ is feasibility cover}{
\Return minimal\textunderscore feasibility\textunderscore cover($Q$)
}
\Else{
\For{$r \in R$}{
    $Q.add(r)$\\
    \If{$Q$ is feasibility cover}{
\Return minimal\textunderscore feasibility\textunderscore cover($Q$)
}
    }
}
\end{algorithm}

\begin{figure}[!b]
  \centering
  \includegraphics[width=\textwidth]{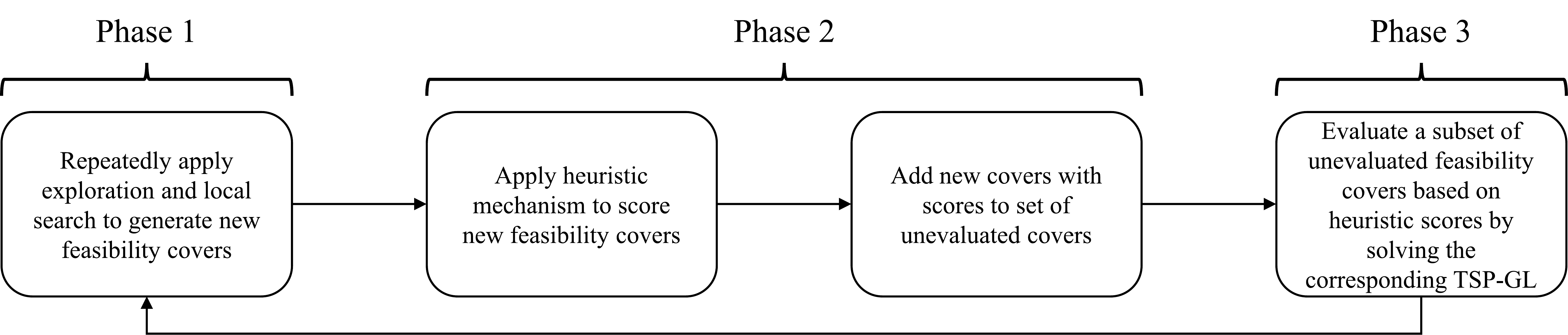}
	\caption{Heuristic flowchart}
	\label{fig:heuristic flowchart}
\end{figure}

In the local search step, the algorithm performs a simple swap procedure on a given feasibility cover. We always perform the local search on two minimal feasibility covers with the smallest upper bounds. Algorithm~\ref{alg: local search} outlines the procedure for one minimal feasibility cover $Q'$. Here, we aim to find a new minimal feasibility cover $Q$ similar to $Q'$. To do so, we assign $Q \gets Q'$ (l. 1) and remove random elements from $Q$ until $Q$ is no longer a feasibility cover (l.3-4). Simultaneously, we add all removed requests to $R$ (l.4). Afterwards, we randomly add new requests that have not been removed before to $Q$ (l.6-11) until either $Q$ is a feasibility cover again (l.8-9) or there are no more requests to add. If $Q$ is a feasibility cover, we determine a minimal feasibility cover through Algorithm~\ref{alg: get_minimum_covers} and return it (l.12-13). Otherwise, we randomly add previously removed requests to $Q$ until it is a feasibility cover again (l.15-22). We then determine a minimal feasibility cover through Algorithm~\ref{alg: get_minimum_covers} and return it.
For all new-found minimal feasibility covers in the exploration and local search step, we use the same heuristic mechanism as in our B\&P approach to score the covers and calculate upper and lower bounds. 
In the solving step, we choose the cover with the lowest upper bound from our queue and solve Problem \ref{TSP-GL ref} on it if its lower bound is smaller than the upper bound of the incumbent. We then update the upper bound of the node by the objective value of the TSP-GL and update the incumbent if necessary.
Figure~\ref{fig:heuristic flowchart} shows the illustrative flowchart for this local search-based heuristic approach. Compared to our B\&P approach, only the initialization through an exploration step and first phase, in which we generate minimal feasibility covers, change. The termination criterion for this approach is a time limit.

\subsection{Hybrid approach}
\label{subsec: hybrid}
Based on our B\&P and heuristic approaches, we can now combine the two of them to develop a hybrid approach. Accordingly, we integrate our heuristic approach's exploration and local search into our B\&P approach. To do so, we initialize our B\&P with a set of minimal feasibility covers $\mathcal{X}_{\theta, \rho}'$ which we determine through our exploration operator. Additionally, after each pricing problem in our B\&P approach, we apply our local search and exploration operator to the new found minimal feasibility cover. This allows us to find minimal feasibility covers similar to the ones determined in our pricing problem through the local search and minimal feasibility covers, which visit a similar number of nodes as the newfound cover through the exploration. Afterwards, we continue with the second and third phases of our B\&P approach. 



\subsection{Speed-up techniques}
\label{sec: Implementation}
We use the following speed-up techniques to enhance the performance of our algorithms. First, closing the optimality gap of the continuous solution of the root node of our B\&P and hybrid approaches is already challenging. Consequently, we limit the number of solved pricing problems that determine new feasibility covers of our algorithm (Figure~\ref{fig:Flowchart SSTSP-GL}) to $n = 5$ in each iteration of the column generation. Here, we also determine an overall lower bound of the STSP-GL in each iteration by calculating the sum of the current objective value of our RMP and the respective objective value of our pricing problem.

Second, after solving $n$ pricing problems, for each feasibility cover $Q$ whose corresponding decision variable $\chi_Q$ is in the basis of the last solution of the RMP, we calculate the upper and lower bound of their corresponding TSP-GL.
We then add each of the selected feasibility covers $Q$ to a priority queue, which sorts them based on their upper bounds. 

Third, in both our B\&P and hybrid approaches, we check whether the feasibility cover's lower bound is smaller than our incumbent's upper bound for the first $m = 5$ feasibility covers in the priority queue. If this holds, we solve the corresponding TSP-GL. Otherwise, we continue to the next element of the priority queue without solving the TSP-GL. We do not test additional feasibility covers in an iteration, even if we solve no corresponding TSP-GLs. 

Fourth, determining the set of all feasibility covers that are not minimal ($\Psi$) for Constraints~\eqref{minimal feasibility covers constraint} is computationally not tractable. Accordingly, we can implement them as lazy constraints. In our computational effectiveness analysis, we did not see a computational advantage of adding these constraints to the pricing problem and sometimes, they even worsened the lower bound of our B\&P and hybrid approaches. Accordingly, we excluded them from the results of our computational analysis.

Fifth, due to the structure of our Benders approach's feasibility cuts~\eqref{feasibility cut benders}, we can cache all feasibility cuts for a set of nodes. If we solve the TSP-GL for a new feasibility cover whose set of nodes was included in another already solved TSP-GL, we initialize our Benders approach with the corresponding feasibility cuts. Finally, we can terminate our Benders approach if the objective value of Problem~\ref{TSP-GL ref} is greater than the current upper bound of our algorithm.

\section{Computational analysis}
\label{sec: Results}
Our analysis is fourfold. First, we analyze the computational efficiency of our algorithms. Second, we study the behavior of upper and lower bounds of our different approaches to understand what drives the complexity of our problem. Third, we analyze the benefits of our stochastic planning approach compared to a deterministic planning approach. Fourth, we aim to understand how the predefined service level and the feasibility frequency impact the design cost. 

To do so, we first compare the bounds of our B\&P, heuristic, and hybrid approach against a MIP benchmark, which solves the deterministic equivalent of our STSP-GL (Problem~\ref{CC-SGMLP decomp}) with a standard solver. Second, we analyze the convergence of bounds of our three approaches against the same MIP benchmark. Third, we compare our stochastic approach against a deterministic approach on multiple instances, varying in graph size, service level, and feasibility frequency. Fourth, we show the difference in design cost on instances with varying service level and feasibility frequency.

\subsection{Instances}
We use the \href{https://w1.cirrelt.ca/~errico/#Instances}{TSPGL2} data set to test our algorithms. The \href{https://w1.cirrelt.ca/~errico/#Instances}{TSPGL2} instance set is an open source data set that builds on the \href{http://comopt.ifi.uni-heidelberg.de/software/TSPLIB95/}{TSPLIB} data set. For analysis (i), we solve instances on complete graphs with 17, 21, 29, 42, 51, 70, and 76 nodes to analyze our algorithms on different graph sizes and analyze our approaches on instances with 20 and 50 scenarios for $\theta = 0.95$ and $\rho = 0.05$ as well as $\theta = 0.9$ and $\rho = 0.1$. We generate five instances considering random seeds for each network dimension, scenario size, and service/feasibility level. Accordingly, this TSPGL2-based instance set consists of 140 different instances. Furthermore, we set $\alpha = 0.25$. For analysis (ii), we focus on one representative example of bound convergence on a graph with 51 nodes for 20 and 50 scenarios. Other graphs from the TSPGL2 data set show similar results to those displayed here. For analysis (iii), we solve instances on complete graphs with 17, 21, 29, 42, 51, and 70 with 20 scenarios for $\theta \in [0.05, 0.1]$ and $\rho \in [0.9, 0.95]$. We again generate five instances considering random seeds for each network dimension, service, and feasibility level. Accordingly, this TSPGL2-based instance set consists of 120 different instances. For analysis (iv), we again focus on one representative example with 51 nodes and 50 scenarios to analyze the impact of service and feasibility level on the design cost. In this context, we solve the STSP-GL for $\theta \in [0.8, 0.85, 0.9, 0.95, 1.0]$ and $\rho \in [0.0, 0.05, 0.1, 0.15, 0.2]$. Other graphs from the TSPGL2 data set again show similar results to those displayed here.

\subsection{Computational setting}
We set the maximum runtime of all algorithms to 60 minutes. In our B\&P and hybrid approach, we limit the computation time for each RMP to 40 minutes and for each pricing problem to 20 minutes. Furthermore, for fair comparison, we terminate the RMP and pricing problem when reaching an optimality gap of two percent. This is necessary, as it may happen that both the RMP and pricing problem cannot find an optimal solution within the maximum computation time. Furthermore, we terminate all algorithms except our heuristic if an optimality gap of two percent is reached.

All our experiments have been conducted on a standard desktop computer equipped with two \texttt{Intel(R) Core(TM) i9-9900, 3.1 GHz CPU} and a total of \texttt{20 GB} of \texttt{RAM}, running \texttt{Ubuntu 20.04}. We have implemented all algorithms in \texttt{Python (3.8.11)} using \texttt{Gurobi 9.5}. Our source code can be found on \href{https://github.com/tumBAIS/STSP-GL}{https://github.com/tumBAIS/STSP-GL}. 

\subsection{Computational effectiveness of proposed methods}
To understand whether the proposed methods are computationally efficient, we benchmark them against a MIP solver solving Problem~\ref{CC-SGMLP decomp}. We do this by aggregating results on instances according to different parameter values. Namely, the number of scenarios, service level, and feasibility frequency.

Tables~\ref{tab: case study 1h} - \ref{tab: case study 1h 0.9 0.1 50 scenarios} show average upper bounds (UB), lower bounds (LB), and optimality gaps (Gap) for the MIP benchmark, our B\&P and hybrid approaches. Additionally, it shows the number of instances for which no primal feasible solution could be found (No UB) for the MIP benchmark and the upper bounds (UB) for our heuristic approach. All data displays the average over five randomly generated instances with a maximum computation time of one hour. We also tested our algorithms with a maximum computation time of two hours, but all our three algorithms and the MIP benchmark rarely improved upper or lower bounds in this additional hour.

Tables~\ref{tab: case study 1h} - \ref{tab: case study 1h 0.9 0.1 50 scenarios} show that for smaller instances with less than 50 nodes, the upper bounds of our B\&P, heuristic, and hybrid approaches are nearly always 0-3\% better than our MIP benchmark. The only outliers are the instances with 42 nodes in Table~\ref{tab: case study 1h 0.9 0.1 50 scenarios}, in which the upper bounds of our B\&P and hybrid approaches are 13\% and of our heuristic are 8\% better. We also see that our B\&P and hybrid approaches find good lower bounds for these instances even though they are nearly always 2-3\% worse than the lower bounds of our MIP benchmark. The only large outliers are the instances with 17 nodes in Table~\ref{tab: case study 1h 0.9 0.1 50 scenarios}, in which the lower bounds of our B\&P and hybrid approaches are 10\% worse than our MIP benchmark. For larger instances with more than 50 nodes, our B\&P and hybrid approaches find much better upper bounds than our MIP benchmark ranging from improvements between 10-44\% in the instances with 51 and 76 nodes and improvements between 74-85\% in instances with 70 nodes.

Furthermore, our B\&P and hybrid approaches find much better lower bounds than our MIP benchmark for the two largest instances with 70 and 76 nodes, with improvements ranging between 8-55\% compared to our MIP benchmark, but cannot close the optimality gap of 5-34\%. Additionally, our MIP benchmark cannot find feasible solutions for at least three of the five problems of our largest instance with 76 nodes in all settings.

If we compare our different approaches over the different instance settings, we see that either our B\&P or our hybrid approach performs better than our heuristic approach on all problems except the ones with 76 nodes, $\theta = 0.95$, $\rho = 0.05$ and 50 scenarios. Furthermore, we see that for most problems, our B\&P and hybrid approaches perform similarly in terms of upper and lower bounds. The main difference in upper and lower bounds is in Table~\ref{tab: case study 1h 0.9 0.1}, where our hybrid approach finds worse lower bounds than our B\&P approach and the fact that our hybrid approach finds better upper bounds than our B\&P approach for the largest instances with 76 nodes in all tables. One explanation for the latter observation is that our hybrid approach finds more minimal feasibility covers it can evaluate before solving the corresponding TSP-GL of one of them. As the instances grow in size, solving the corresponding TSP-GL becomes computationally expensive, i.e., determining the upper and lower bounds of more minimal feasibility covers leads to a higher probability of choosing a good minimal feasibility cover to evaluate.

\begin{table}[!h]
  \centering
  \caption{TSPGL2-based instances with 20 scenarios $\theta = 0.95$, $\rho = 0.05$}
  \label{tab: case study 1h}
  \small
\begin{tabular}{ccccccccccccc}
\multicolumn{2}{l}{Instance} &  \multicolumn{4}{c}{MIP} &  \multicolumn{3}{c}{B\&P} &  \multicolumn{3}{c}{Hybrid} & Heuristic \\
\cmidrule(r){1-2} \cmidrule(lr){3-6} \cmidrule(lr){7-9} \cmidrule(lr){10-12} \cmidrule(l){13-13}
 Nodes  & $\vert S \vert$ & UB & LB & Gap [\%] & No UB & UB & LB & Gap [\%] & UB & LB & Gap [\%] & UB\\
\midrule
       17 &        20 &            3231 &            3187 &               1.38 &            0 &            3230 &            3230 &               0.00 &               3230 &               3230 &                  0.00 &                  3230 \\
       21 &        20 &            4258 &            4239 &               0.45 &            0 &            4253 &            4215 &               0.91 &               4252 &               4218 &                  0.82 &                  4252 \\
       29 &        20 &            2531 &            2489 &               1.67 &            0 &            2524 &            2412 &               4.67 &               2524 &               2426 &                  4.05 &                  2523 \\
       42 &        20 &            1046 &             985 &               6.29 &            0 &            1036 &             990 &               4.99 &               1036 &                975 &                  6.42 &                  1050 \\
       51 &        20 &             668 &             566 &              18.28 &            0 &             591 &             560 &               5.59 &                591 &                560 &                  5.65 &                   613 \\
       70 &        20 &            5086 &             686 &             637.77 &            0 &             925 &             885 &               4.89 &                927 &                821 &                 12.95 &                   988 \\
       76 &        20 &             864 &             561 &              51.18 &            3 &             720 &             671 &               7.35 &                718 &                667 &                  7.70 &                   778 \\
\bottomrule
\end{tabular}
\end{table}
\begin{table}[!h]
  \centering
  \caption{TSPGL2-based instances with 50 scenarios $\theta = 0.95$, $\rho = 0.05$}
  \label{tab: case study 1h 50 scenarios}
  \small
\begin{tabular}{ccccccccccccc}
\multicolumn{2}{l}{Instance} &  \multicolumn{4}{c}{MIP} &  \multicolumn{3}{c}{B\&P} &  \multicolumn{3}{c}{Hybrid} & Heuristic \\
\cmidrule(r){1-2} \cmidrule(lr){3-6} \cmidrule(lr){7-9} \cmidrule(lr){10-12} \cmidrule(l){13-13}
 Nodes  & $\vert S \vert$ & UB & LB & Gap [\%] & No UB & UB & LB & Gap [\%] & UB & LB & Gap [\%] & UB\\
\midrule
       17 &        50 &            3279 &            3224 &               1.72 &            0 &            3276 &            3169 &               3.46 &               3276 &               3276 &                  0.00 &                  3276 \\
       21 &        50 &            4259 &            4224 &               0.83 &            0 &            4257 &            4205 &               1.24 &               4255 &               4225 &                  0.73 &                  4255 \\
       29 &        50 &            2536 &            2495 &               1.64 &            0 &            2530 &            2407 &               5.11 &               2526 &               2399 &                  5.33 &                  2530 \\
       42 &        50 &            1079 &            1024 &               5.42 &            0 &            1075 &             964 &              11.58 &               1071 &                978 &                  9.55 &                  1073 \\
       51 &        50 &             694 &             596 &              16.36 &            0 &             615 &             598 &               2.72 &                615 &                600 &                  2.60 &                   629 \\
       70 &        50 &            6079 &             692 &             779.12 &            0 &             962 &             885 &               8.90 &                966 &                866 &                 11.57 &                   992 \\
       76 &        50 &             933 &             451 &              63.18 &            4 &             828 &             702 &              18.06 &                784 &                697 &                 12.55 &                   778 \\
\bottomrule
\end{tabular}
\end{table}
\begin{table}[!h]
  \centering
  \caption{TSPGL2-based instances with 20 scenarios $\theta = 0.9$, $\rho = 0.1$}
  \label{tab: case study 1h 0.9 0.1}
  \small
\begin{tabular}{ccccccccccccc}
\multicolumn{2}{l}{Instance} &  \multicolumn{4}{c}{MIP} &  \multicolumn{3}{c}{B\&P} &  \multicolumn{3}{c}{Hybrid} & Heuristic \\
\cmidrule(r){1-2} \cmidrule(lr){3-6} \cmidrule(lr){7-9} \cmidrule(lr){10-12} \cmidrule(l){13-13}
 Nodes  & $\vert S \vert$ & UB & LB & Gap [\%] & No UB & UB & LB & Gap [\%] & UB & LB & Gap [\%] & UB\\
\midrule
       17 &        20 &            3132 &            3096 &               1.16 &            0 &            3130 &            2896 &               8.31 &               3131 &               2955 &                  6.34 &                  3203 \\
       21 &        20 &            3945 &            3904 &               1.03 &            0 &            3937 &            3706 &               6.33 &               3937 &               3763 &                  4.78 &                  3987 \\
       29 &        20 &            2166 &            2129 &               1.75 &            0 &            2167 &            2068 &               4.85 &               2167 &               2047 &                  5.86 &                  2238 \\
       42 &        20 &             983 &             808 &              21.74 &            0 &             947 &             947 &               0.00 &                945 &                845 &                 12.74 &                   981 \\
       51 &        20 &             582 &             456 &              27.83 &            0 &             527 &             484 &               9.37 &                520 &                467 &                 11.66 &                   581 \\
       70 &        20 &            3154 &             516 &             511.39 &            0 &             827 &             671 &              24.71 &                836 &                630 &                 32.71 &                   953 \\
       76 &        20 &            1123 &             403 &             171.74 &            3 &             631 &             576 &              10.04 &                630 &                539 &                 16.88 &                   738 \\
\bottomrule
\end{tabular}
\end{table}
\begin{table}[!h]
  \centering
  \caption{TSPGL2-based instances with 50 scenarios $\theta = 0.9$, $\rho = 0.1$}
  \label{tab: case study 1h 0.9 0.1 50 scenarios}
  \small
\begin{tabular}{ccccccccccccc}
\multicolumn{2}{l}{Instance} &  \multicolumn{4}{c}{MIP} &  \multicolumn{3}{c}{B\&P} &  \multicolumn{3}{c}{Hybrid} & Heuristic \\
\cmidrule(r){1-2} \cmidrule(lr){3-6} \cmidrule(lr){7-9} \cmidrule(lr){10-12} \cmidrule(l){13-13}
 Nodes  & $\vert S \vert$ & UB & LB & Gap [\%] & No UB & UB & LB & Gap [\%] & UB & LB & Gap [\%] & UB\\
\midrule
       17 &        50 &            3269 &            3223 &               1.43 &            0 &            3268 &            2892 &              13.03 &               3268 &               2913 &                 12.18 &                  3277 \\
       21 &        50 &            3967 &            3905 &               1.59 &            0 &            3963 &            3786 &               4.80 &               3964 &               3792 &                  4.64 &                  4007 \\
       29 &        50 &            2191 &            2152 &               1.82 &            0 &            2200 &            2098 &               4.93 &               2203 &               2092 &                  5.35 &                  2293 \\
       42 &        50 &            1109 &             806 &              36.75 &            1 &             962 &             962 &               0.00 &                962 &                904 &                  6.97 &                  1014 \\
       51 &        50 &             689 &             471 &              46.55 &            0 &             552 &             506 &               9.33 &                554 &                492 &                 12.56 &                   595 \\
       70 &        50 &            6081 &             322 &          243816.66 &            0 &             886 &             672 &              31.68 &                896 &                671 &                 33.49 &                   955 \\
       76 &        50 &             832 &             255 &              97.54 &            4 &             785 &             571 &              37.80 &                742 &                567 &                 30.97 &                   744 \\
\bottomrule
\end{tabular}
\end{table}

Furthermore, we see that the number of scenarios as well as the reduction of the service level and feasibility frequency drives complexity, i.e., leads to worse optimality gaps in our B\&P and hybrid approaches. The reason for that is, that with a growing number of scenarios as well as shrinking service level and feasibility frequency, the number of minimal feasibility covers grows significantly and more TSP-GLs have to be solved.

\begin{figure}[!t]
\vspace{-2ex}
	\begin{subfigure}[t]{.3\textwidth}
		\centering		
		\includegraphics[width=\textwidth]{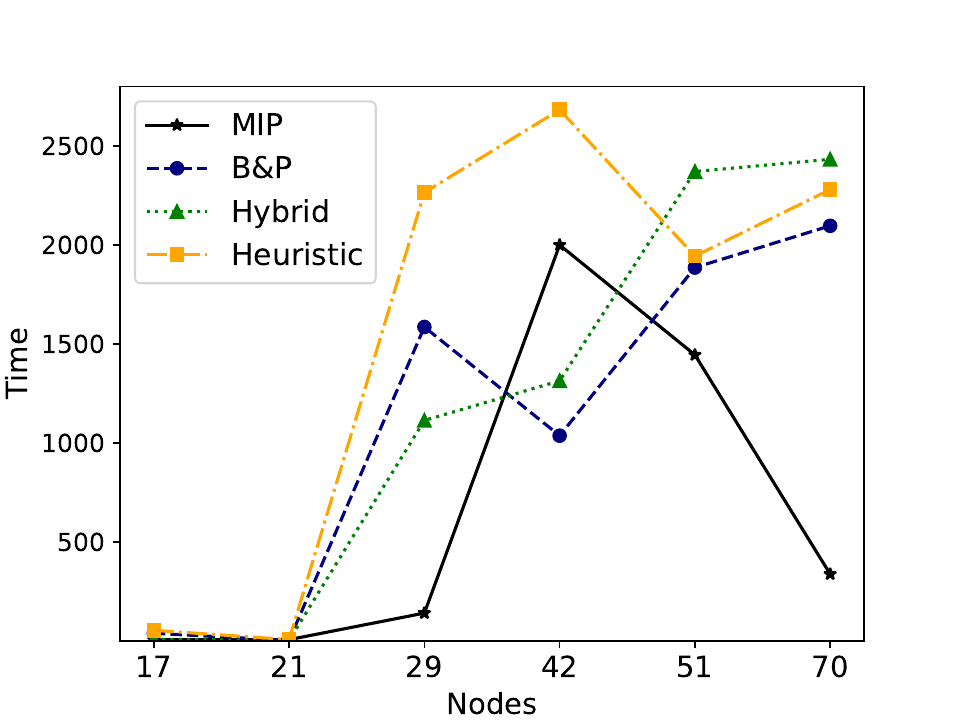}
        \captionsetup{format=hang}
		\caption{Computation times until final upper bound}
		\label{fig: runtime time}
	\end{subfigure}
	\begin{subfigure}[t]{.3\textwidth}
		\centering
		\includegraphics[width=\textwidth]{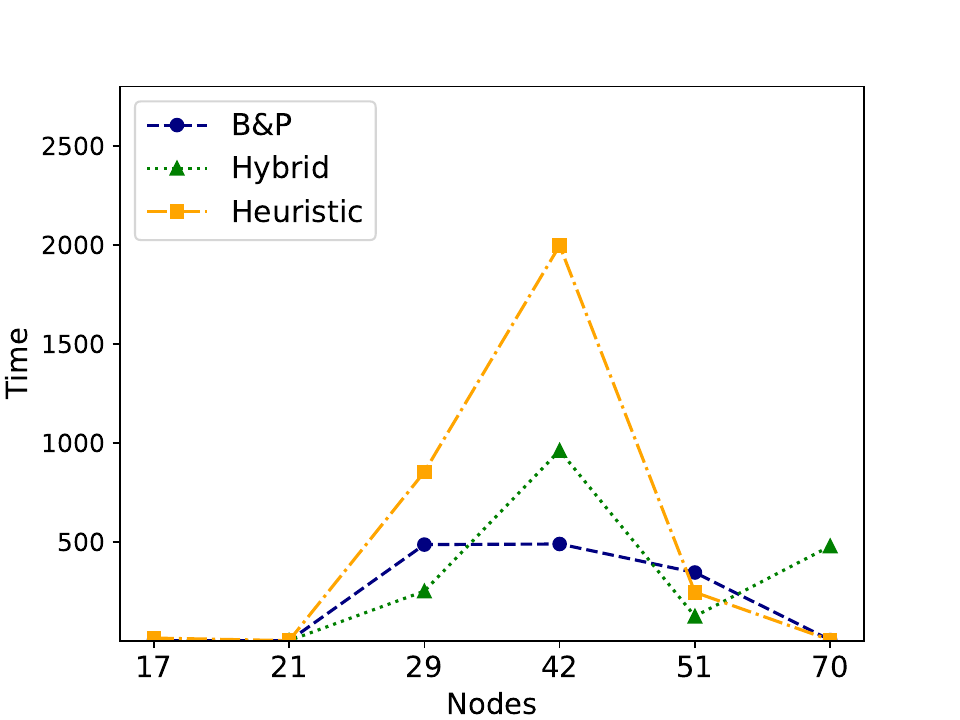}
        \captionsetup{format=hang}
		\caption{Computation times until better upper bound than best MIP solution}
		\label{fig: runtime MIP time}
	\end{subfigure}
 	\begin{subfigure}[t]{.3\textwidth}
		\centering
		\includegraphics[width=\textwidth]{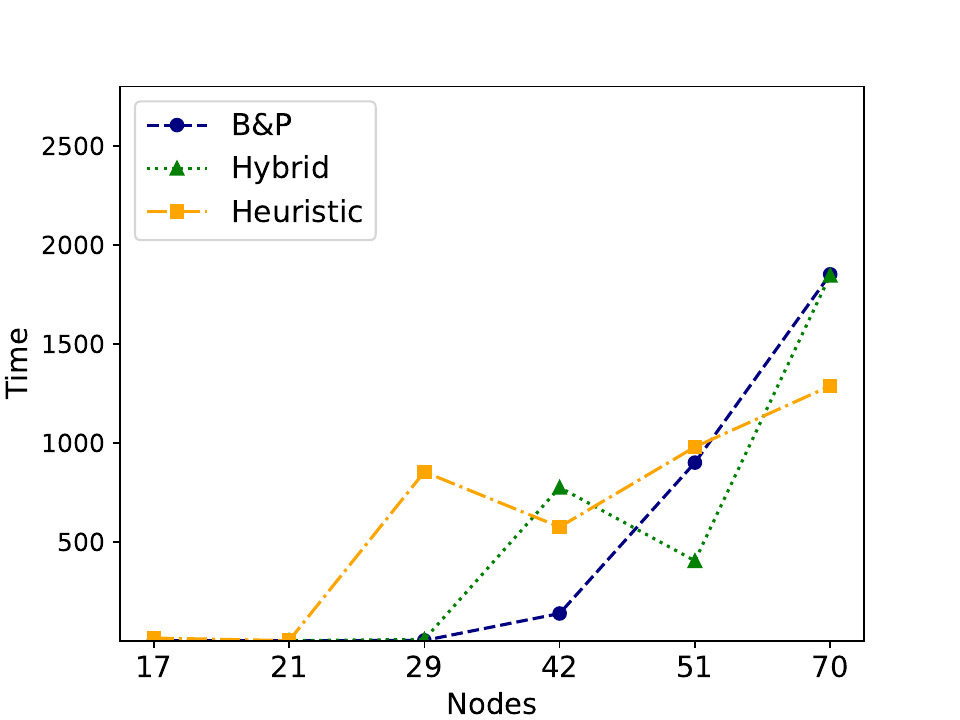}
        \captionsetup{format=hang}
		\caption{Computation times until upper bound that is not more than 1\% worse than final upper bound}
		\label{fig: runtime 2 percent}
	\end{subfigure}
	\caption{Computation times of TSPGL2-based instances with 20 scenarios $\theta = 0.95$, $\rho = 0.05$}
	\label{fig: case study 1h 20 scenarios computation times}
\end{figure}

Figure~\ref{fig: case study 1h 20 scenarios computation times} shows average computation times until the final upper bound is found (a), average computation times until a better upper bound than the best MIP bound is found in our B\&P, hybrid, and heuristic approaches (b), and average computation times until an upper bound is found that is not more than 1\% worse than the final upper bound (c) for five instances with 20 scenarios with $\theta = 0.95$ and $\rho = 0.05$. The results for instances with 50 scenarios or $\theta = 0.9$ and $\rho = 0.1$ are similar to those shown here. We do not add the results for instances with 76 nodes, as the MIP does not find solutions for most of them.
We see that our B\&P approach, on average, finds better upper bounds than our MIP in under 500 seconds for all instances. Additionally, we see that there may be a significant difference between finding an upper bound that is within one percent of the final upper bound and finding the final upper bound, e.g., instances with 29 nodes with 1585 seconds (Figure~\ref{fig: runtime time}) versus 2 seconds (Figure~\ref{fig: runtime 2 percent}). We observe similar behavior for our hybrid and heuristic approach. Additionally, we see that except for the largest instances with 70 nodes, either our B\&P or our hybrid approaches find better upper bounds than the MIP faster than our heuristic approach. 

\subsection{Convergence of bounds}
\begin{figure}[!t]
\centering
\vspace{-2ex}
	\begin{subfigure}[t]{.45\textwidth}
		\centering		
		\includegraphics[width=\textwidth]{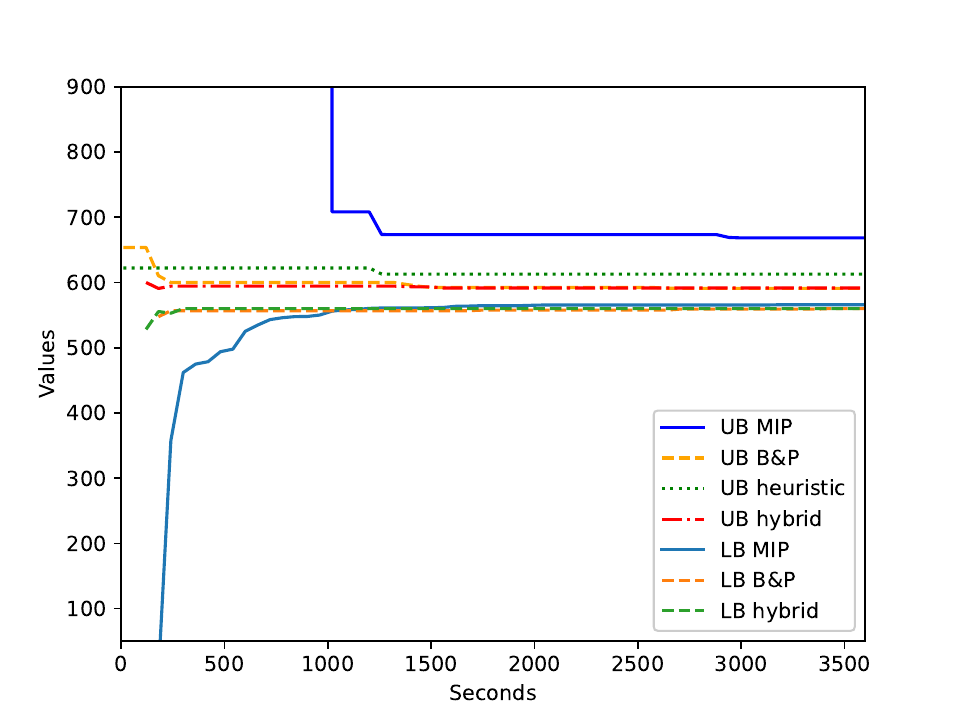}
        \captionsetup{format=hang}
		\caption{Convergence of upper and lower bound for 51 nodes and 20 scenarios}
		\label{fig: bound convergence small}
	\end{subfigure}
	\begin{subfigure}[t]{.45\textwidth}
		\centering
		\includegraphics[width=\textwidth]{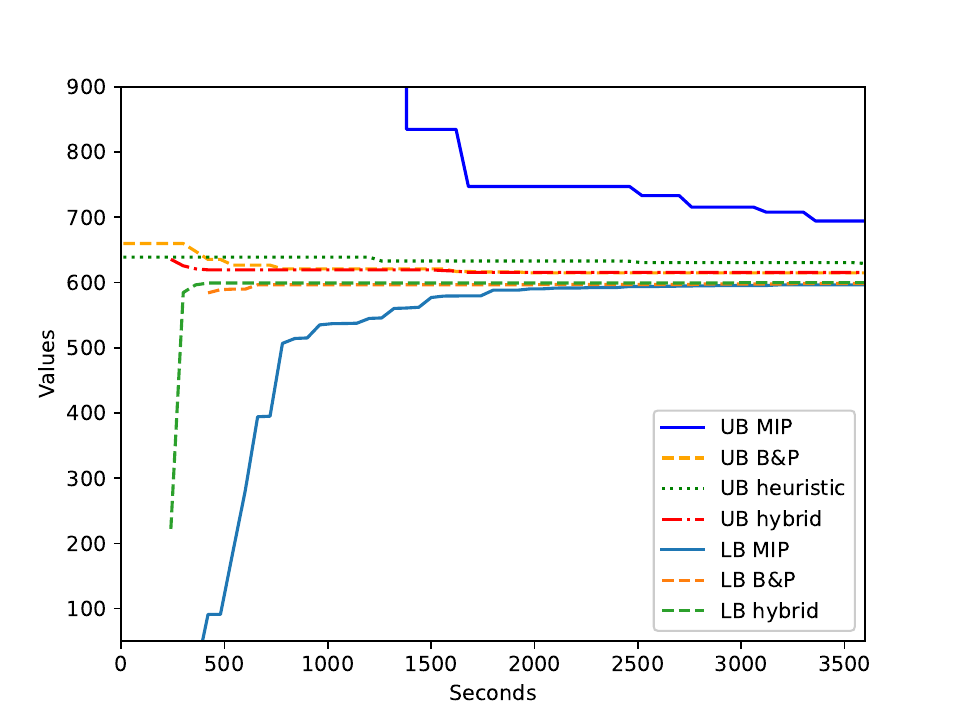}
        \captionsetup{format=hang}
		\caption{Convergence of upper and lower bound for\\ 51 nodes and 50 scenarios}
		\label{fig: bound convergence large}
	\end{subfigure}
	\caption{Average bound convergence in instances with 20 and 50 scenarios for $\theta = 0.95$, $\rho = 0.05$}
	\label{fig: bound convergence 0.95}
\end{figure}

In this section, we analyze the convergence of the upper and lower bounds of our B\&P, heuristic, and hybrid approaches with our MIP benchmark on the TSPGL2 dataset. Note that our heuristic has no lower bound.
Figure~\ref{fig: bound convergence 0.95} shows the average convergence of upper and lower bounds for instances with 51 nodes and a service level of $\theta = 0.95$ and $\rho = 0.05$. We see that our B\&P and hybrid approaches find better upper bounds faster than our MIP benchmark for instances with 20 scenarios over the entire computation time and quickly find better upper bounds than our heuristic approach. Our heuristic approach also finds better upper bounds than our MIP benchmark. We also see that all methods ultimately converge to similar lower bounds.
For instances with 50 scenarios, we again find upper bounds with our B\&P and hybrid approaches much faster than our MIP benchmark. Our heuristic approach finds good upper bounds even quicker than our B\&P and hybrid approaches but cannot improve significantly over the computation time. Furthermore, we see that with more scenarios, our B\&P and hybrid approaches find upper bounds much faster than our MIP benchmark but cannot close the optimality gap.

\begin{figure}[!t]
\centering
\vspace{-2ex}
	\begin{subfigure}{.45\textwidth}
		\centering		
		\includegraphics[scale=0.5]{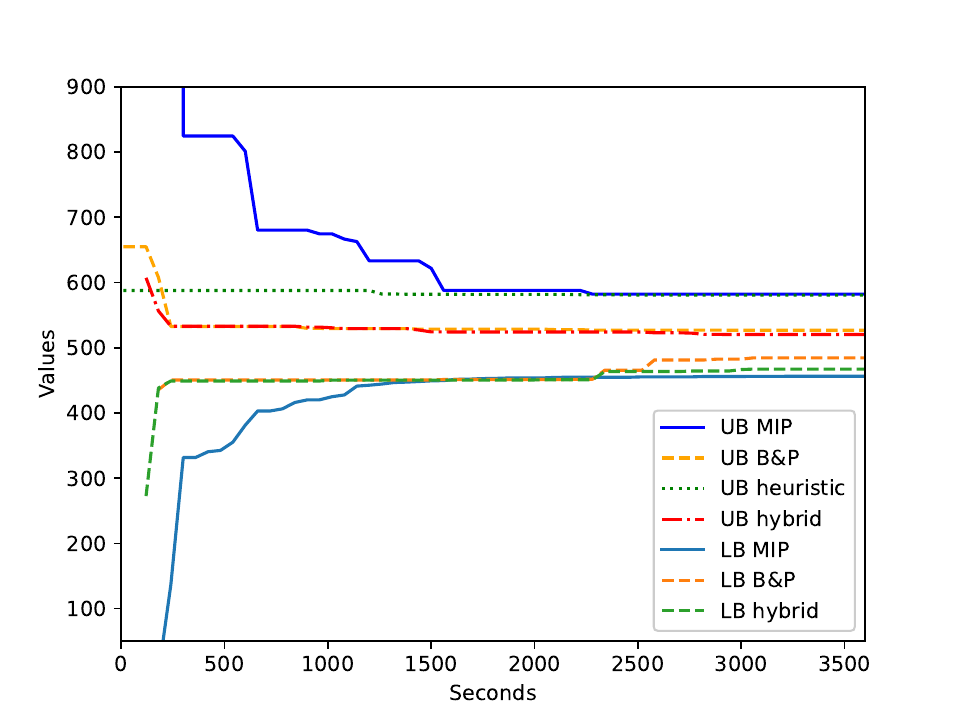}
        \captionsetup{format=hang}
		\caption{Convergence of upper and lower bound for 51 nodes and 20 scenarios}
		\label{fig: bound convergence small}
	\end{subfigure}
	\begin{subfigure}{.45\textwidth}
		\centering
		\includegraphics[scale=0.5]{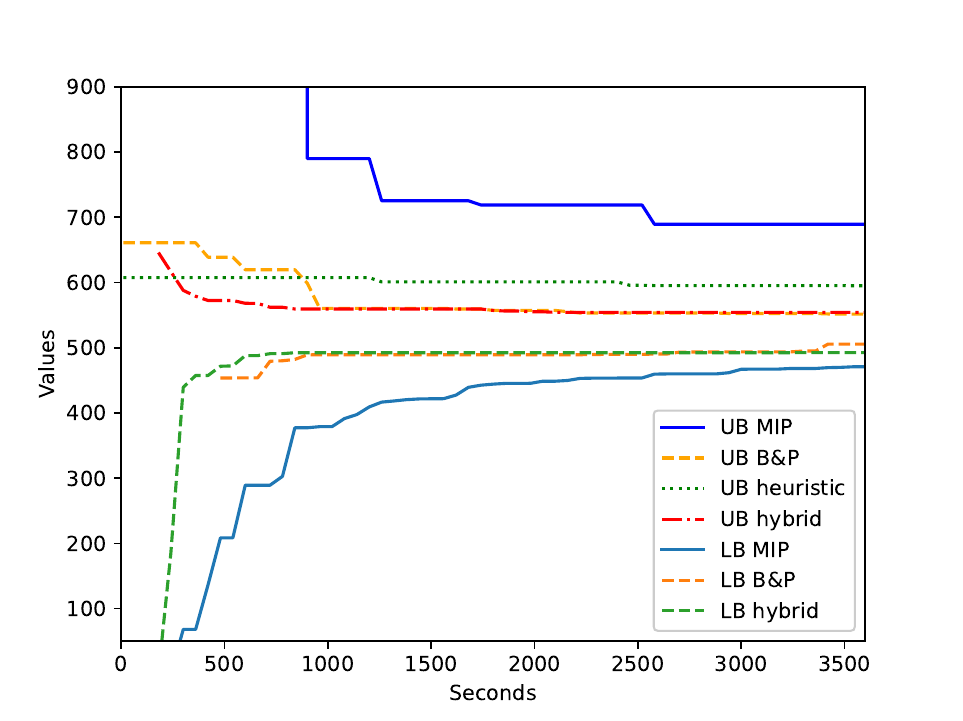}
        \captionsetup{format=hang}
		\caption{Convergence of upper and lower bound for 51 nodes and 50 scenarios}
		\label{fig: bound convergence large}
	\end{subfigure}
	\caption{Average bound convergence instances with 20 and 50 scenarios for $\theta = 0.9$, $\rho = 0.1$}
	\label{fig: bound convergence 0.9}
\end{figure}

Figure~\ref{fig: bound convergence 0.9} shows the average convergence of upper and lower bounds for instances with 51 nodes and a service level of $\theta = 0.9$ and $\rho = 0.1$. We see that for 20 scenarios, our B\&P and hybrid approaches find better upper bounds faster than our MIP benchmark. Similar to the higher service level, our heuristic approach finds upper bounds fast but cannot improve much during the rest of the run. Still, contrary to the higher service level instances with 20 scenarios, our MIP benchmark cannot find better bounds than our heuristic approach. Looking at the lower bounds of our B\&P and hybrid approaches compared to our MIP benchmark, we see that for a lower service level, our B\&P and hybrid approaches find equally good lower bounds to our MIP benchmark. 
For instances with 50 scenarios, our MIP benchmark takes much longer than our B\&P, heuristic, and hybrid approaches to find feasible solutions. Additionally, our B\&P and hybrid approaches find better lower bounds faster than our MIP benchmark. Again, none of the algorithms can close the optimality gap within one hour.

With this analysis of bound convergence, we see that our B\&P, heuristic, and hybrid approaches are more robust against an increase in scenarios than our MIP benchmark and can find feasible solutions of good quality much faster.


\subsection{Value of Stochastic Solution}
It is common for solutions determined by deterministic models to derive better objective values than solutions determined by stochastic models when we evaluate them in a deterministic setting. At the same time, the deterministic solutions are often infeasible or derive worse objective values than stochastic solutions when we evaluate them in a stochastic setting. In this section, we show that this property also holds for the STSP-GL.

\begin{table}[!b]
\centering
\caption{Comparison between deterministic and stochastic approach}
\label{tab: comparison deterministic and stochastic approach}
\small
\begin{tabular}{cccccccccccc}
\multicolumn{3}{l}{Instance} &  \multicolumn{5}{c}{Deterministic approach} &  \multicolumn{4}{c}{Stochastic approach}\\
\cmidrule(r){1-3} \cmidrule(lr){4-8} \cmidrule(lr){9-12}
 Nodes &  $\theta$ &  $\rho$ &  Design cost & $\bar{N}$ &  $\bar{d}$ &  $\bar{\rho}$ &  Inf. &  Design cost & $\bar{N}$ &  $\bar{d}$ & $\bar{\rho}$ \\
\midrule
       17 &  0.95 & 0.05 &                  2846.1 &              0.91 &                    0.96 &                          0.31 &                           5 / 5 &      3083.7 &  0.99 &        0.97 &              0.05 \\
       21 &  0.95 & 0.05 &                  3714.3 &              0.90 &                    0.95 &                          0.43 &                           5 / 5 &      4060.5 &  1.00 &        0.98 &              0.05 \\
       29 &  0.95 & 0.05 &                  2077.5 &              0.90 &                    0.95 &                          0.45 &                           5 / 5 &      2464.8 &  0.99 &        0.98 &              0.04 \\
       42 &  0.95 & 0.05 &                   937.8 &              0.90 &                    0.95 &                          0.41 &                           5 / 5 &       989.4 &  0.93 &        0.95 &              0.05 \\
       51 &  0.95 & 0.05 &                   543.6 &              0.85 &                    0.95 &                          0.42 &                           5 / 5 &       585.3 &  0.90 &        0.95 &              0.05 \\
       70 &  0.95 & 0.05 &                   888.0 &              0.89 &                    0.95 &                          0.47 &                           5 / 5 &       916.2 &  0.92 &        0.95 &              0.05 \\
\midrule       
       17 &  0.95 & 0.10 &                  2846.1 &              0.91 &                    0.96 &                          0.31 &                           5 / 5 &      3083.7 &  0.99 &        0.94 &              0.10 \\
       21 &  0.95 & 0.10 &                  3714.3 &              0.90 &                    0.95 &                          0.43 &                           5 / 5 &      4060.5 &  1.00 &        0.96 &              0.10 \\
       29 &  0.95 & 0.10 &                  2077.5 &              0.90 &                    0.95 &                          0.45 &                           5 / 5 &      2365.5 &  0.97 &        0.96 &              0.10 \\
       42 &  0.95 & 0.10 &                   937.8 &              0.90 &                    0.95 &                          0.41 &                           5 / 5 &       983.1 &  0.92 &        0.94 &              0.10 \\
       51 &  0.95 & 0.10 &                   540.6 &              0.85 &                    0.95 &                          0.42 &                           5 / 5 &       573.3 &  0.89 &        0.94 &              0.10 \\
       70 &  0.95 & 0.10 &                   888.0 &              0.89 &                    0.95 &                          0.47 &                           5 / 5 &       911.4 &  0.91 &        0.93 &              0.10 \\
\midrule       
       17 &  0.90 & 0.05 &                  2611.8 &              0.85 &                    0.91 &                          0.42 &                           5 / 5 &      3018.3 &  0.96 &        0.95 &              0.05 \\
       21 &  0.90 & 0.05 &                  3380.4 &              0.86 &                    0.90 &                          0.41 &                           5 / 5 &      3908.1 &  0.96 &        0.93 &              0.05 \\
       29 &  0.90 & 0.05 &                  1812.3 &              0.83 &                    0.90 &                          0.44 &                           5 / 5 &      2114.1 &  0.90 &        0.91 &              0.05 \\
       42 &  0.90 & 0.05 &                   876.0 &              0.80 &                    0.90 &                          0.43 &                           5 / 5 &       924.0 &  0.85 &        0.91 &              0.05 \\
       51 &  0.90 & 0.05 &                   485.7 &              0.76 &                    0.90 &                          0.43 &                           5 / 5 &       520.8 &  0.80 &        0.90 &              0.05 \\
       70 &  0.90 & 0.05 &                   797.4 &              0.77 &                    0.90 &                          0.44 &                           5 / 5 &       828.0 &  0.81 &        0.90 &              0.05 \\
\midrule       
       17 &  0.90 & 0.10 &                  2611.8 &              0.85 &                    0.91 &                          0.42 &                           5 / 5 &      2999.1 &  0.95 &        0.92 &              0.10 \\
       21 &  0.90 & 0.10 &                  3380.4 &              0.86 &                    0.90 &                          0.41 &                           5 / 5 &      3762.0 &  0.92 &        0.91 &              0.10 \\
       29 &  0.90 & 0.10 &                  1812.3 &              0.83 &                    0.90 &                          0.44 &                           5 / 5 &      2069.4 &  0.87 &        0.89 &              0.10 \\
       42 &  0.90 & 0.10 &                   876.0 &              0.80 &                    0.90 &                          0.43 &                           5 / 5 &       904.5 &  0.84 &        0.88 &              0.10 \\
       51 &  0.90 & 0.10 &                   486.9 &              0.76 &                    0.90 &                          0.43 &                           5 / 5 &       512.4 &  0.78 &        0.89 &              0.10 \\
       70 &  0.90 & 0.10 &                   792.6 &              0.77 &                    0.90 &                          0.41 &                           5 / 5 &       814.5 &  0.80 &        0.88 &              0.10 \\
\bottomrule
\end{tabular}
\end{table}

Table \ref{tab: comparison deterministic and stochastic approach} shows the comparison between a deterministic approach and our stochastic approach. Specifically, we compare our stochastic approach to the deterministic setting in which we solve the STSP-GL for only one scenario in which we calculate the mean demand for each request over all original scenarios. In this context, we test the two different approaches on different instances with $\theta \in [0.95, 0.9]$, $\rho \in [0.05, 0.1]$, and 20 scenarios. We solve five instances for each individual configuration. Table \ref{tab: comparison deterministic and stochastic approach} shows the design cost, percentage of nodes in the tour ($\bar{N}$), average demand of served demand ($\bar{d}$), average percentage of scenarios in which the service level is not fulfilled ($\bar{\rho}$), and the number of deterministic solutions which are infeasible in our stochastic setting (Inf.). 

We see that the deterministic approach finds tours with much lower design costs than our stochastic approach as, on average, the tours include fewer nodes. At the same time, as is expected for deterministic approaches in a stochastic setting, we see that for all instances, the solutions of the deterministic approach can not fulfill the service level for more than 30 percent of all scenarios. This makes all deterministic solutions infeasible for our stochastic setting, as for them to be feasible, $\bar{\rho}$ would need to be smaller than $\rho$, i.e., 0.05 and 0.1. Furthermore, we see that the percentage of served demand ($\bar{d}$) is often smaller in the deterministic solutions compared to our stochastic approaches. In this context, it is interesting that $\bar{d}$ is significantly higher than $\theta * (1 - \rho)$ which we expected in the stochastic approach. This phenomenon suggests that to fulfill our feasibility frequency, we have to serve more demand than required by our service level. Furthermore, except for the instances with 29 nodes, $\theta = 0.95$ and $\rho = 0.05$, $\bar{\rho}$ is as expected equal to $(1 - \rho)$ for all stochastic solutions.

Summarizing, our stochastic approach accounts for uncertain passenger demand by including more nodes in the tour, increasing the design cost. Additionally, significantly more demand is served in our stochastic solutions compared to the deterministic solutions.

\subsection{Analysis of service and feasibility level impact on design cost}
Our model presumes a service level, i.e., the market share served, and the frequency with which it is served. To understand which drives costs, we aggregate results on instances according to different parameter values for the service level $\theta$ and the frequency of meeting it $(1 - \rho)$.

\begin{figure}[!b]
\centering
\begin{subfigure}[t]{.45\textwidth}
	\centering
	\includegraphics[width=1.1\textwidth]{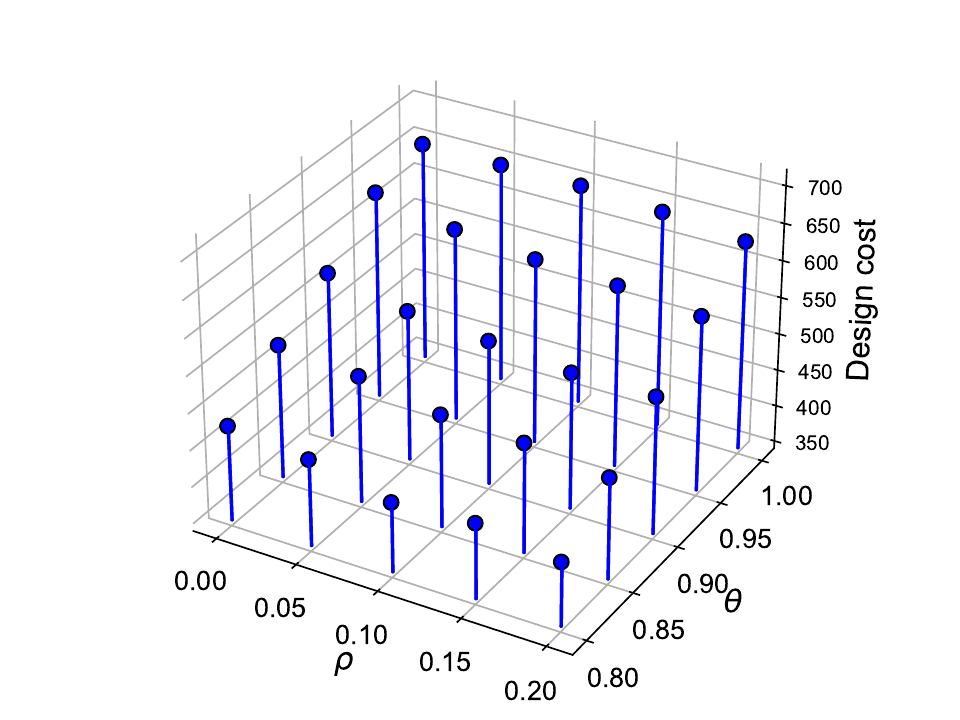}
    \captionsetup{format=hang}
	\caption{Average design cost}
		\label{fig: theta rho impact design}
\end{subfigure}
\begin{subfigure}[t]{.45\textwidth}
	\centering
	\includegraphics[width=1.1\textwidth]{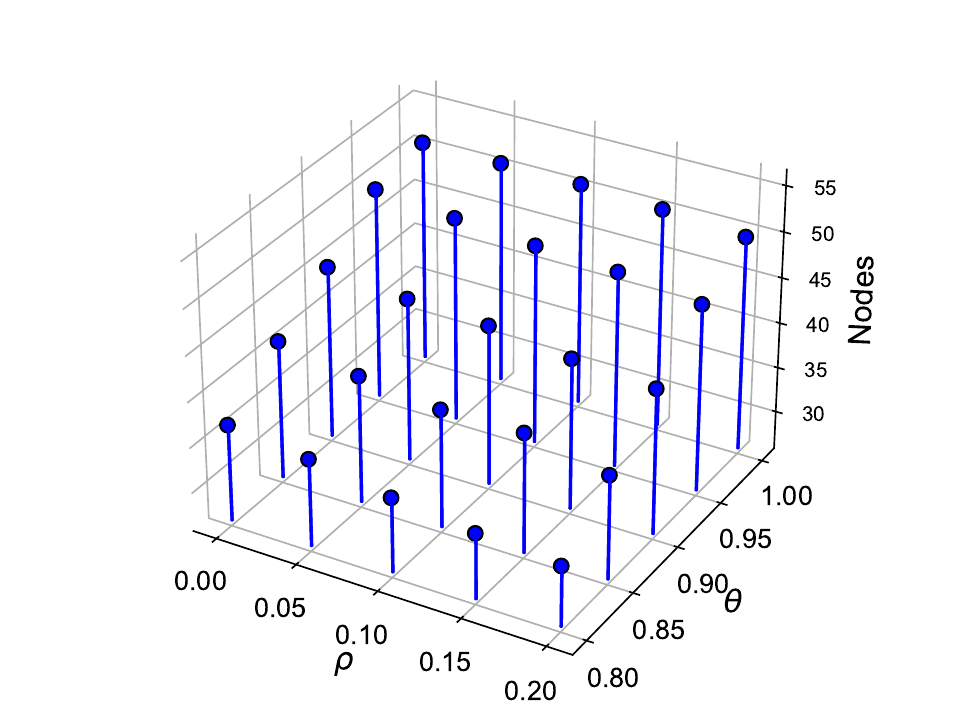}
    \captionsetup{format=hang}
	\caption{Average number of nodes}
		\label{fig: theta rho impact nodes}
\end{subfigure}
\caption{Average design cost and number of nodes for varying $\theta$ and $\rho$ values on instance with 51 nodes}
\label{fig: theta rho impact}
\end{figure}

Figure~\ref{fig: theta rho impact} shows the difference in design costs and number of nodes in the tour for service levels $\theta \in [0.8, 0.85, 0.9, 0.95, 1.0]$ and feasibility levels $\rho \in [0.0, 0.05, 0.1, 0.15, 0.2]$ on the TSPGL2 instance with 51 nodes with 50 scenarios. We see that the design cost of a route depends more on the desired market size $\theta$ than on the frequency of meeting it ($1 - \rho$), i.e., reducing $\theta$ by $0.05$ on average decreases design costs by 8.5\% while reducing $1 - \rho$ by $0.05$ on average only reduces design costs by 2\%. We see a similar reduction in the number of nodes in the tour. Explicitly, reducing $\theta$ by $0.05$ on average decreases the number of nodes in the tour by 9\% while reducing $1 - \rho$ by $0.05$ on average only reduces the number of nodes in the tour by 2\%. This is not surprising as the number of nodes in a tour highly correlates to the design costs of the tour.
Consequently, a DAS operator aiming to reduce the design costs in the tactical planning problem should focus on the desired service level instead of on the frequency of meeting it.

\section{Conclusion and future work}
\label{sec: Conclusion}
With this work, we introduced the STSP-GL, which can be used to find a general ordering of stops in the tactical planning of DAS under uncertain passenger demand. Furthermore, we proposed an algorithmic framework to determine near-optimal solutions. To do so, we modeled the STSP-GL as a stochastic program and determined its deterministic equivalent. This allowed us to utilize a B\&P, hybrid, and heuristic approach, which decompose the STSP-GL into multiple TSP-GLs. We performed extensive computational experimentation to analyze how our different algorithms perform on different graph sizes, number of scenarios, service, and feasibility levels. For this purpose, we generated a large set of instances from the \href{https://w1.cirrelt.ca/~errico/#Instances}{TSPGL2} dataset. We show that our B\&P and hybrid approaches find better upper bounds than our MIP benchmark within 500 seconds for nearly all instance types but struggle with closing the optimality gap. Especially for larger graph sizes, our B\&P and hybrid approaches improve the MIP benchmark's upper bounds by up to 85\%. Furthermore, we showed that our heuristic approach finds feasible solutions within a couple of seconds but is not able to significantly improve on the initial upper bounds during the remainder of the computation time. Additionally, we saw that the additional local search and exploration of our hybrid approach makes it more robust against higher numbers of scenarios compared to our B\&P approach. Furthermore, we showed that even though the design costs of solutions derived through a deterministic approach are lower compared to solutions of a stochastic approach, they are not feasible in a stochastic setting, i.e., the stochastic approach includes more nodes in the tours to meet the demand with a higher probability. Finally, we showed that the impact of the service level on the design cost is higher than the impact of the feasibility level. Here, we are able to, on average, reduce the design costs by 8.5\% when lowering the service level by 5\% while only reducing the design costs by 2\% when lowering the feasibility frequency by 5\%.

Further research may expand on the current methodology of our heuristic to develop a more sophisticated hybrid approach to outperform our B\&P and hybrid approaches. Additionally, our current distribution of requests is nearly uniform, i.e., the probability of request A having positive demand in a scenario is similar to request B having positive demand in a scenario. Investigating the impact of different request distributions on the design and routing costs remains an interesting avenue for future research.

\section*{Acknowledgment}
This work was supported by the German Federal Ministry of Education and Research under Grant 03ZU1105FA.


%
\singlespacing{
\bibliographystyle{model5-names}
\bibliography{ms}} 
\newpage
\onehalfspacing
\begin{appendices}
	\normalsize
\end{appendices}
\end{document}